\documentclass[authoryear,12pt]{elsarticle}
\usepackage{cancel}
\usepackage{caption}
\usepackage{color}
\usepackage{lscape}
\usepackage{afterpage}
\usepackage{pstricks}
\usepackage{pst-plot}
\usepackage{longtable}
\usepackage{dcolumn}
\usepackage{pst-node}
\usepackage{amsmath}
\usepackage{amssymb}
\usepackage{amsthm}
\usepackage{multirow}
\usepackage{colortab}
\usepackage{color}
\usepackage{array}
\usepackage{slashbox}
\usepackage{colortbl}
\usepackage{subfigure}
\usepackage{textcomp}
\usepackage{pstricks, pst-node}
\usepackage{pst-all}
\usepackage{algorithm,algorithmic}
\usepackage{url}
\usepackage{soul}
\usepackage{hyperref}
\psset{arrows=->, labelsep=3pt, mnode=circle}

\usepackage{tabularx}

\usepackage{eurosym}

\newfont{\rams}{msbm10 scaled\magstep1}
\newcommand{\rea}{\mbox{\rams \symbol{'122}}}

\newenvironment{resumeT}{\begin{list}{}{\setlength{\rightmargin}{\leftmargin}}\item[]
{\centering {\bf \it~~~}
\par}\item[]\ignorespaces}{\unskip\end{list}}

\newtheorem{example}{Example}[section]
\begin{document}
\title{Robust stochastic sorting with interacting criteria hierarchically structured}

\author[Eco]{\rm Sally Giuseppe Arcidiacono}
\ead{s.arcidiacono@unict.it}
\author[Eco]{\rm Salvatore Corrente}
\ead{salvatore.corrente@unict.it}
\author[Eco,por]{\rm Salvatore Greco}
\ead{salgreco@unict.it}

\address[Eco]{Department of Economics and Business, University of Catania, Corso Italia, 55, 95129  Catania, Italy}
\address[por]{University of Portsmouth, Portsmouth Business School, Centre of Operations Research and Logistics (CORL), Richmond Building, Portland Street, Portsmouth PO1 3DE, United Kingdom}

\date{}
\maketitle

\vspace{-1cm}

\begin{resumeT}

\textbf{Abstract:} \noindent In this paper we propose a new multiple criteria decision aiding method to deal with sorting problems in which alternatives are evaluated on criteria structured in a hierarchical way and presenting interactions. The underlying preference model of the proposed method is the Choquet integral, while the hierarchical structure of the criteria is taken into account by applying the Multiple Criteria Hierarchy Process. Considering the Choquet integral based on a 2-additive capacity, the paper presents a procedure to find all the minimal sets of pairs of interacting criteria representing the preference information provided by the Decision Maker (DM). Robustness concerns are also taken into account by applying the Robust Ordinal Regression and the Stochastic Multicriteria Acceptability Analysis. Even if in different ways, both of them provide recommendations on the hierarchical sorting problem at hand by exploring the whole set of capacities compatible with the preferences provided by the DM avoiding to take into account only one of them. The applicability of the considered method to real world problems is demonstrated by means of an example regarding rating of European Countries by considering  economic and financial data provided by Standard \& Poor's Global Inc.

{\bf Keywords}: {Decision Support System, Sorting Problems, Interactions between criteria, Parsimonious models, Robustness concerns}
\end{resumeT}

 \pagenumbering{arabic}

\section{Introduction}
Sorting is one of the different type of problems that can be dealt by using a Multiple Criteria Decision Aiding (MCDA) method \citep{GrecoEhrgottFigueira2016}. In sorting problems a set of alternatives have to be assigned to a set of classes ordered from the worst to the best with respect to the preferences provided by the Decision Maker (DM) (see \citealt{ZopounidisDoumpos2002} for a survey on MCDA sorting methods). Several sorting methods have been developed to cope with such type of problems having at the basis different preference models such as value functions \citep{KeeneyRaiffa1976}, outranking relations \citep{BransVincke1985,FigueiraEtAl2013} or Decision Rules \citep{GrecoMatarazzoSlowinski2001}. For example, without any ambition to be exhaustive, let us list the following sorting methods mentioning the underlying preference model: AHPSort \citep{IshizakaPearmanNemery2012} and UTADIS \citep{DevaudGroussaudJacquetLagreze1980,ZopounidisDoumpos1999} are based on an additive value function; ELECTRE Tri B \citep{Yu1992}, ELECTRE Tri C \citep{AlmeidaDiasFigueiraRoy2010} and ELECTRE Tri-nC \citep{AlmeidaDiasFigueiraRoy2012} are based on the ELECTRE methods, while FlowSort \citep{NemeryLamboray2008} is based on PROMETHEE methods; finally, \citep{GrecoMatarazzoSlowinski2002} deal with a sorting problem using the Dominance Based Rough Set Approach (DRSA). 

While all these methods can be distinguished between compensatory and non-compensatory ones, depending on the fact that they take into account compensation or non-compensation between criteria, none of them can deal with another important characteristic of decision making problems being the possible interaction between criteria. In simple terms, two criteria are positively interacting if the importance assigned to them together is greater than the sum of the importance assigned to each of them separately, while two criteria are negatively interacting if the importance assigned to these criteria together is lower than the sum of the importance assigned to them singularly. In this paper, we will therefore develop a multicriteria sorting method based on the Choquet integral preference model \citep{Choquet1953}, that is, the most well-known non-additive integral used in literature to take into account such a type of interactions (see \citealt{Grabisch1996} for the application of the Choquet integral in MCDA). To the best of our knowledge, very few contributions are proposing methods to deal with sorting problems with interacting criteria and a comparison between our approach and these methods will be recalled in Section \ref{ComparisonSection} after that the new proposal will be described in detail. 

The application of the different sorting methods mentioned above, as well as of all MCDA methods, involves the knowledge of several parameters. For example, on one hand, the application of an additive value function implies the knowledge of tradeoffs between criteria or, in general, the shape of the marginal value functions; on the other hand, outranking relations can be computed by knowing the weights of criteria as well as some thresholds used to take into account the uncertainty, ambiguity or imprecision of the alternatives' evaluations \citep{RoyFigueiraAlmeidaDias2014}. Such parameters can be obtained by asking the DM to provide a direct or an indirect preference information, depending on the fact that he is willing and able to provide an exact value for them or to provide some preference examples from which parameters compatible with these preference can be inferred. In general, the indirect way of providing preference information is preferable for the DM since it involves a lower cognitive effort \citep{JacquetLagrezeSiskos2001}. Sorting methods do not represent an exception with respect to this aspect and, consequently, different contributions take into account preference information provided by the DM in an indirect way to discover one or more set of instances of the considered preference model compatible with them. For example, different types of preference information provided by the DM in an indirect way can be found in \cite{CorrenteEtAl2017a,GrecoMousseauSlowinski2010,KadzinskiEtAl2016,KadzinskiCiomekSlowinski2015,KoksalanOzpeynirci2009} for sorting methods based on value functions and in \cite{DiasMousseau2003,KadzinskiCiomek2016,KadzinskiSlowinski2013,KadzinskiTervonenFigueira2015,LeroyMousseauPirlot2011,ZhangEtAl2014} for sorting methods based, instead, on outranking relations. \\
In this paper, the proposed sorting method based on the Choquet integral preference model will take into account an indirect preference information provided by the DM. To build a model being able to represent the preferences provided by the DM but, at the same time, as parsimonious as possible \citep{ArcidiaconoCorrenteGreco2020}, we introduce a procedure to define a minimal sets of pairs of interacting criteria and for the core of these sets. We think that this is an important contribution of the paper since the core is composed of the pairs of interactions being really necessary to explain the information provided by the DM. Moreover, the idea of the minimal sets of pairs of interacting criteria can be used in all the applications of the Choquet integral, not only those ones related to sorting problems.

Two more aspects of decision making problems are taken into account and introduced in the proposed method, that is, the hierarchy of criteria and the robustness concerns. 

Many sorting problems present alternatives evaluated on criteria structured in a hierarchical way. For such a reason, in addition to the comprehensive classification of the alternatives considering all criteria together, a partial classification of the alternatives at hand on the basis of set of homogeneous criteria constituting macrocriteria can be beneficial for the DM who has the possibility to get a finer representation of the goodness of the same alternatives. The hierarchy of criteria will be taken into account in the proposed method by the Multiple Criteria Hierarchy Process \citep{CorrenteGrecoSlowinski2012} which permits the DM to provide preferences and to get information not only at the global level but also considering single aspects of the decision problem he is dealing with. 

Considering an indirect preference information provided by the DM, more than one instance of the preference model can be compatible with them. Therefore, giving an assignment of the alternatives at hand using only one of them can be reductive or, anyway, misleading. To overcome this drawback, the proposed method provides assignments by exploring the whole set of instances of the assumed preference model by applying the Robust Ordinal Regression (ROR; \citealt{GrecoMousseauSlowinski2008}) and the Stochastic Multicriteria Acceptability Analysis (SMAA; \citealt{LahdelmaHokkanenSalminen1998}). On one hand, the application of the ROR to the proposed sorting method will produce necessary and possible assignments of the considered alternatives, being assignments that hold for all or for at least one compatible model; on the other hand, the application of SMAA will give back information in terms of frequencies with which an alternative will be assigned to a particular class or to an interval of possible classes. Of course, the integration of ROR and SMAA with the MCHP will permit to the DM to get the illustrated information not only in a comprehensive way and, therefore, considering all criteria in a while, but also considering subsets of criteria of particular interest for him. Moreover, to get a final classification on the basis of the frequencies provided by the SMAA application, a new procedure based on the minimization of misclassification errors is presented for the first time in this paper. 

The applicability of the introduced hierarchical, interacting and robust sorting method will be demonstrated by means of an example regarding the financial rating of European sovereign debts evaluated by Standard \& Poor’s Global Inc. In particular, we shall show that the proposed method is quite predictive with respect to other sorting methods based on different preference models. 

The paper is structured in the following way. In the next section, a first didactic example showing the necessity to take into account interacting effects between criteria in sorting problems is given; Section \ref{prel} provides the background of the new proposal. In particular, the hierarchical Choquet integral preference model is briefly recalled; in Section \ref{SortingMethod}, we introduce the new hierarchical, interacting and robust multicriteria sorting method; a comparison between the proposed method and others taking into account the interactions between criteria is given in Section \ref{ComparisonSection}; Section \ref{CaseStudy} contains the case study in which the new method is applied to the financial rating of 28 European Countries; finally, conclusions and further directions of research are collected in Section \ref{concl}.

\section{A first motivating example}

Let suppose that a Credit Rating Agency (CRA)\footnote{A credit rating agency is a company that assigns a rating to a debtor according to its ability to pay back its debt.} has to provided an ordinal classification of several sovereign bonds evaluated on economic, governmental and financial points of view. The bonds have to be classified in four classes ordered from the best to the worst as follows: $AA\succ A\succ BB\succ B,$ where $X \succ Y$ means that the alternatives from class X are preferred to the alternatives from class Y, $X,Y=AA,A,BB,B$. From an MCDA perspective this is a sorting problem where the bonds are the alternatives to be classified and the economic, governmental and financial points of view represent the evaluation criteria ($G=\{ {g}_{Eco},{g}_{Gov},{g}_{Fin} \}$). The CRA is therefore interested into building a sorting model that could be used to classify all the considered bonds. 

For this reason, the analyst decides to apply the simplest preference model being an additive value function $U:A\rightarrow[0,1]$

\begin{equation}\label{AdditiveFunction}
U(x) = U(g_{1}(x),\ldots,g_{n}(x))=\displaystyle\sum_{i=1}^{n}u_{i}(g_{i}(x))
\end{equation}

\noindent where $A$ denotes the set of alternatives, $g_{1},\ldots,g_{n}$ are the considered criteria and $u_i$ are the marginal value functions related to criteria $g_i$. Without loss of the generality, in the following we suppose that $g_i:A \rightarrow \rea_+$. We shall assume that the greater $U(x)$ the better the alternative $a$ and that classes are ordered in an increasing way with respect to preferences so that $C_{1}$ and $C_{p}$ denote, respectively, the worst and the best classes to which an alternative can be assigned. Of course if $a$ is assigned to $C_{h}$ and $b$ is assigned to $C_{k}$, with $h>k$, then $U(a)>U(b)$. \\
To build such a model, let us suppose that the analyst asks to the representative of the CRA to assign four bonds $\{a,b,c,d\}$ evaluated on a [0,20] scale on the three aspects above as shown in Table \ref{exampevaltab}. Consequently, the representative of the CRA performs the assignments shown in the last column of the same table. 

\begin{table}[h]
\begin{center}
\caption{Evaluation of the four sovereign bonds and the CRA assignment\label{exampevaltab}}	
		\begin{tabular}{c|ccc|c}
		\hline
		\hline
		Sovereign bonds$\backslash$criteria & $g_{Eco}$  & $g_{Gov}$  & $g_{Fin}$ & Class \\
		\hline
		\hline
		$a$      &  11    &    9    &   5   &  $BB$ \\
		$b$      &   7    &   12    &   5   &  $B$ \\
		$c$      &  11    &    9    &   8   &  $A$ \\
		$d$      &   7    &   12    &   8   &  $AA$ \\
		\hline
		\hline
		\end{tabular}
\end{center}
\end{table}

Since $a$ is assigned to $BB$ and $b$ is assigned to $B$, following what has been said before, $U(a)>U(b)$. Analogously, because of the assignments of $c$ and $d$ to $A$ and $AA$, respectively, it should also hold that $U(d)>U(c)$. Using the function in eq. (\ref{AdditiveFunction}), the two pieces of information are therefore translated into the following constraints 
$$\left\{
\begin{array}{l}
{U(a)={u}_{Eco}(11) + {u}_{Gov}(9) + {u}_{Fin}(5) > {u}_{Eco}(7) + {u}_{Gov}(12) + {u}_{Fin}(5)}=U(b); \\[2mm]
{U(c)={u}_{Eco}(11) + {u}_{Gov}(9) + {u}_{Fin}(8)} < {u}_{Eco}(7) + {u}_{Gov}(12) + {u}_{Fin}(8)=U(d),
\end{array}
\right.
$$
\noindent being clearly in contradiction.\\
Looking carefully at the performances shown in Table \ref{exampevaltab}, it is indeed evident that the set of criteria $\left\{g_{Eco},g_{Gov}\right\}$ is not preferentially independent of criterion $g_{Fin}$ \citep{Wakker1989}. Indeed, $a$ and $c$ have the same performances on $g_{Eco}$ and $g_{Gov}$ (11 and 9) as well as $b$ and $d$ (7 and 12). At the same time, $a$ and $b$ have the same performance on $g_{Fin}$ (5) as happens also for $c$ and $d$ (8). Consequently, if $\left\{g_{Eco},g_{Gov}\right\}$ would be preferentially independent of $g_{Fin}$, the preference of $a$ over $b$ should imply the preference of $c$ over $d$ even if, the provided assignments of $c$ and $d$ imply exactly the opposite preference. This is due to the fact that there is a certain degree of interaction between the criteria at hand. 

The very didactic example shown above proves that eventual positive or negative interactions between criteria need to be properly taken into account in MCDA methods developed to deal with sorting problems.

\section{Background}\label{prel}
In this section, we shall briefly introduce the notation that will be used along the paper, together with the methods that constitute the basis of the new proposal, that is, the Multiple Criteria Hierarchy Process (MCHP; \citealt{CorrenteGrecoSlowinski2012}) and the Choquet integral preference model \citep{Choquet1953,Grabisch1996}. In the MCHP, all criteria are not at the same level but they are structured in a hierarchical way. The application of the MCHP will give the possibility to consider a particular aspect of the considered problem, without taking into account all criteria simultaneously. The Choquet integral, instead, is able to take into consideration the possible positive or negative interactions between criteria. The description of the Choquet integral will be given in the context of the MCHP.   

\subsection{Notation}
The used notation is the following:
\begin{itemize}
\item $A=\left\{a,b,\ldots\right\}$ is the set of alternatives,
\item ${G}$ is the whole set of criteria (at all levels of the hierarchy), and $g_\mathbf{0}$ is the root criterion, 
\item ${I}_{{G}}$ is the set of indices of the criteria in $G$,
\item ${G}_{{E}}\subseteq{G}$ is the set of all elementary criteria in ${G}$, that is, the criteria located at the bottom of the hierarchy. The alternatives will be evaluated on these criteria only. Consequently, to each $a\in A$ will be associated the vector $\left(g_{\mathbf{t}_1}(a),\ldots,g_{\mathbf{t}_{|{G}_{E}|}}(a)\right)\in\rea^{|{G}_{E}|},$ which components are the evaluations of $a$ on the considered elementary criteria. In the following, without loss of generality, we shall assume that all elementary criteria have an increasing direction of preference that is, the greater the evaluation of an alternative on the considered elementary criterion, the better the alternative is on it (for some recent contributions taking into account non-monotonic criteria in sorting problems, see \cite{GuoLiaoLiu2019,KadzinskiEtAl2020,LiuEtAl2019}), 
\item ${E}_{G}\subseteq{I}_{G}$ is the set of indices of the elementary criteria,
\item $g_{\mathbf{r}}$ is a non-elementary criterion, that is a criterion $g_{\mathbf{r}}\in{G}\setminus{G}_{{E}}$. In particular, by $g_{\mathbf{0}}$ we denote the root criterion and, consequently, the whole set of criteria in the hierarchy,
\item ${E}(g_{\mathbf{r}})\subseteq{E}_{G}$ is the set of indices of all the elementary criteria descending from $g_{\mathbf{r}}$ (it follows that ${E}(g_{\mathbf{0}})={E}_{G}$),
\item given ${F}\subseteq{G}$, ${E}({F})=\displaystyle\cup_{g_{\mathbf{r}}\in{F}}{E}(g_{\mathbf{r}})$, that is the set of the indices of the elementary criteria descending from at least one criterion in $F$,
\item ${G}^{k}_{\mathbf{r}}\subseteq{G}$ is the set of subcriteria of $g_{\mathbf{r}}$ located at level $k$.
\end{itemize}

\subsection{The hierarchical Choquet integral preference model}\label{MCHPChoquet}%
As already explained in the introduction, in real world applications, the set of evaluation criteria is not mutually preferentially independent since the criteria at hand can present a certain degree of positive or negative interaction. In such cases, to aggregate the performances of the alternatives on the criteria at hand, non additive integrals are used and, among them, the Choquet integral is the most well-known. In the following, we shall briefly describe its application considering the set of elementary criteria $G_{E}$. 

The Choquet integral application is based on a capacity $\mu: 2^{{G}_{{E}}}\rightarrow[0,1]$, being a set function such that $\mu(\emptyset)=0,\; \mu({G}_{E})=1$ and $\mu(R)\leqslant \mu(T),\;\forall R \subseteq T \subseteq G_{E}$. Given $a\in A$, the Choquet integral of its performances' vector $\left(g_{\mathbf{t}_1}(a),\ldots,g_{\mathbf{t}_{|{G}_{E}|}}(a)\right)$ is
\begin{equation*}
Ch\left(g_{\mathbf{t}_1}(a),\ldots,g_{\mathbf{t}_{|{G}_{E}|}}(a)\right)=Ch(a)=\overset{|G_{E}|}{\underset{i=1}{\sum }} \mu \left( N_{\mathbf{t}_i}\right) \left[ g_{(\mathbf{t}_i)}(a) - g_{(\mathbf{t}_{i-1})}(a) \right],
\end{equation*}

\noindent where $_{(\cdot)}$ stands for a permutation of the indices of the elementary criteria such that $0=g_{\left(\mathbf{t}_{0}\right)}(a)\leqslant g_{\left(\mathbf{t}_{1}\right)}(a) \leqslant\ldots\leqslant  g_{\left(\mathbf{t}_{|G_{E}|}\right)}\left(a\right)$ and $N_{\mathbf{t}_i}=\left\{g_{\mathbf{t}_j}\in G_{E}: g_{\mathbf{t}_j}(a)\geqslant g_{(\mathbf{t}_i)}(a)\right\}.$ In the following, for the sake of simplicity, we shall say ``the Choquet integral of $a$" instead of ``the Choquet integral of $\left(g_{\mathbf{t}_1}(a),\ldots,g_{\mathbf{t}_{|{G}_{E}|}}(a)\right)$". 

To make things easier, in general, the M\"{o}bius transform $m$ of the capacity $\mu$ is taken into account \citep{Rota}. $m:2^{G_{E}}\rightarrow \rea$ is a set function such that, for all $R\subseteq G_{E}$, $\mu (R)=\displaystyle\sum_{T \subseteq R}m(T)$ and, vice versa, $\displaystyle m(R)=\sum_{T\subseteq R}(-1)^{(|R|-|T|)}\mu(T).$

In terms of the M\"{o}bius transform $m$ of $\mu$, the Choquet integral of $a$ can be computed as

\begin{equation}\label{MobiusChoquet}
Ch(a)=\displaystyle\sum_{T\subseteq G_{E}}m(T)\min_{g_{\mathbf{t}_i}\in T}\{g_{\mathbf{t}_i}(a)\},     
\end{equation}

\noindent while the monotonicity and normalization constraints of $\mu$ can be rewritten as 

\begin{itemize}
\item $m(\emptyset)=0,$ $\displaystyle\sum_{T\subseteq G_E}m(T)=1$,
\item $\forall g_{\mathbf{t}_i}\in G_{E}$ and $\forall R\subseteq G_{E}\setminus\{g_{\mathbf{t}_i}\}$, $\displaystyle\sum_{T\subseteq R}m\left(T\cup\{g_{\mathbf{t}_i}\}\right)\geqslant 0$.
\end{itemize}

Since the application of the Choquet integral involves the knowledge of $2^{|G_{E}|}-2$ parameters, that is, one for each subset of criteria of $G_{E}$ apart from $\emptyset$ and $G_{E}$ ($\mu(\emptyset)=0$ and $\mu(G_{E})=1$), in real world applications, in general, $2$-additive measures are used where, a measure $\mu$ is said $k$-additive \citep{grabisch1997k} if its M\"{o}bius transform is such that $m(T)=0$ for all $T\subseteq G_{E}: \;|T|>k$. By using a $2$-additive measure, the Choquet integral can be reformulated as 

\begin{equation}\label{MinsChoquet}
Ch(a)=\displaystyle\sum_{g_{\mathbf{t}_i}\in G_{E}}m\left(\{g_{\mathbf{t}_i}\}\right)g_{\mathbf{t}_i}(a)+\sum_{\{g_{\mathbf{t}_i},g_{\mathbf{t}_j}\}\subseteq G_{E}}m\left(\{g_{\mathbf{t}_i},g_{\mathbf{t}_j}\}\right)\min\{g_{\mathbf{t}_i}(a),g_{\mathbf{t}_j}(a)\}
\end{equation}

\noindent while the monotonicity and normalization constraints become 

\begin{itemize}
\item $m(\emptyset)=0$, $\displaystyle\sum_{g_{\mathbf{t}_i}\in G_E}m\left(\{g_{\mathbf{t}_i}\}\right)+\sum_{\{g_{\mathbf{t}_i},g_{\mathbf{t}_j}\}\subseteq G_{E}}m\left(\{g_{\mathbf{t}_i},g_{\mathbf{t}_j}\}\right)=1,$
\item $\forall g_{\mathbf{t}_i}\in G_{E}$ and $\forall T\subseteq G_{E}\setminus\{g_{\mathbf{t}_i}\},$ $\;m\left(\{g_{\mathbf{t}_i}\}\right)+\displaystyle\sum_{g_{\mathbf{t}_j}\in T}m\left(\{g_{\mathbf{t}_i},g_{\mathbf{t}_j}\}\right)\geqslant 0.$
\end{itemize}

 Applying the MCHP to the Choquet integral preference model, for each non-elementary criterion $g_{\mathbf{r}}$ it is possible to define a capacity $\mu_{\mathbf{r}}^{k}$ (in the following, we shall write only $\mu_{\mathbf{r}}$) on the power set of $G_{\mathbf{r}}^{k}$ that assigns a value to each subset of criteria descending from $g_{\mathbf{r}}$ and located at the level $k$ \citep{angilella2016MCHP}. The capacity $\mu_{r}$ and its M\"{o}bius transform $m_{r}$ can be expressed in terms of $\mu$ and $m$ defined on the power set of $G_{E}$ (see \citealt{angilella2016MCHP} for technical details) and, consequently, also the Choquet integral of $a$ on $g_{\mathbf{r}}$ can be expressed in terms of $\mu$ by applying the following equality

\begin{equation}
Ch_{\textbf{r}}(a)=\frac{Ch({a}_{\textbf{r}})}{\mu({E}(g_{\mathbf{r}}))}
\label{ChoquetHierarchy}
\end{equation}

\noindent where ${a}_{\textbf{r}}$ is a fictitious alternative having the same evaluations of $a$ on all elementary criteria descending from $g_{\mathbf{r}}$ (that is criteria $g_{\mathbf{t}}$ with $\mathbf{t}\in {E}(g_{\mathbf{r}})$) and null evaluation on the remaining elementary criteria.

Using the Choquet integral preference model, it is obvious that the importance of a criterion is not dependent on itself only, but also from its contribution to all possible coalitions of criteria. For such a reason, the Shapley index \citep{shapley} and the Murofushi index \citep{Murofushi1993} are computed obtaining, on one hand, the importance of a criterion and, on the other hand, the importance of a pair of criteria.  \\
Generalizing these two indices to the MCHP, given a criterion $g_{\left(\mathbf{r},w\right)}\in G_{\mathbf{r}}^{k}$, that is a subcriterion of $g_{\mathbf{r}}$ at the level $k$, the Shapley index $\varphi_{\mathbf{r}}\left(\left\{g_{\left(\mathbf{r},w\right)}\right\}\right)$ is computed as 

\begin{small}
\begin{equation}\label{shapley_MCHP_2}
\varphi_{\mathbf{r}}\left(\left\{g_{\left(\mathbf{r},w\right)}\right\}\right)=\left(\;\sum_{\mathbf{t}\in E\left(g_{\left(\mathbf{r},w\right)}\right)}m\left(\{g_{\mathbf{t}}\}\right)+\sum_{\mathbf{t_1},\mathbf{t_2}\in E\left(g_{\left(\mathbf{r},w\right)}\right)}m\left(\{g_{\mathbf{t_1}},g_{\mathbf{t_2}}\}\right)+\sum_{\substack{\mathbf{t_{1}}\in E\left(g_{\left(\mathbf{r},w\right)}\right) \\ \mathbf{t_{2}}\in E\left({G}_{\mathbf{r}}^{k}\setminus\{g_{\left(\mathbf{r},w\right)}\}\right)}}\frac{m(\{g_{\mathbf{t_{1}}},g_{\mathbf{t_{2}}}\})}{2}\right)\frac{1}{\mu(E(g_{\mathbf{r}}))},
\end{equation}
\end{small}
\noindent and, analogously, given $g_{\left(\mathbf{r},w_1\right)},g_{\left(\mathbf{r},w_2\right)}\in G_{\mathbf{r}}^{k}$, the Murofushi index $\varphi_{\mathbf{r}}\left(\{g_{\left(\mathbf{r},w_{1}\right)},g_{\left(\mathbf{r},w_{2}\right)}\}\right)$ is computed as 

$$
\varphi_{\mathbf{r}}\left(\{g_{\left(\mathbf{r},w_{1}\right)},g_{\left(\mathbf{r},w_{2}\right)}\}\right)=\left(\displaystyle\sum_{\substack{\mathbf{t_1}\in E\left(g_{\left(\mathbf{r},w_1\right)}\right)\\\mathbf{t_2}\in E\left(g_{\left(\mathbf{r},w_2\right)}\right)}}m\left(\{g_{\mathbf{t_{1}}},g_{\mathbf{t_{2}}}\}\right)\right)\frac{1}{\mu\left(E\left(g_{\mathbf{r}}\right)\right)}.
$$

\section{A hierarchical, interacting and robust sorting method}\label{SortingMethod}
In multiple criteria decision problems which are hierarchically structured, for each non-elementary criterion $g_{\mathbf{r}}$, the sorting procedure consists in  the assignment of each alternative $a\in A$ to one class $C_h$, $h\in\left\{1,\ldots,p_{\mathbf{r}}\right\}$, taking into account only the elementary criteria descending from $g_\mathbf{r}$, where $C_{p_{\mathbf{r}}}$ is the class of top performing alternatives and $C_1$ is the class of the worst alternatives. Of course, the number of class to which an alternative can be assigned can be dependent on the considered criterion and has not to be the same for all criteria.   

According to \cite{CorrenteEtAl2017a}, we consider a \textit{threshold based} sorting procedure. This means that each class $C_h$ is defined by a lower threshold $b^{\mathbf{r}}_{h-1}$ and an upper threshold $b^{\mathbf{r}}_{h}$. To perform the assignment, we need, therefore, $p_{\mathbf{r}}+1$ thresholds $b_{0}^{\mathbf{r}},\ldots,b_{p_{\mathbf{r}}}^{\mathbf{r}}$ for each criterion $\mathbf{r}$, such that $0=b^{\mathbf{r}}_{0}<\ldots<b^{\mathbf{r}}_{p_{\mathbf{r}}}$ and $b^{\mathbf{r}}_{p_{\mathbf{r}}}$ is the highest possible evaluation of an alternative on $g_{\mathbf{r}}$. This means that $b_{p_{\mathbf{r}}}^{\mathbf{r}}=Ch_{\mathbf{r}}\left(opt_{\mathbf{t}_1},\ldots,opt_{\mathbf{t_r}}\right)$, where $\mathbf{t}_{1},\ldots,\mathbf{t}_{r}\in E\left(g_{\mathbf{r}}\right)$ are the indices of the elementary criteria descending from $g_{\mathbf{r}}$ and $opt_{\mathbf{t}}=\displaystyle\max_{a\in A}g_{\mathbf{t}}(a)$ for all $\mathbf{t}\in\{\mathbf{t}_{1},\ldots,\mathbf{t}_{r}\}$.

For each alternative $a\in A$, for each non-elementary criterion $g_{\mathbf{r}}$ and for each class $C_h$, $h=1,\ldots,p_{\mathbf{r}}$:
\begin{itemize}
\item $a$ is assigned to $C_{h}$ on $g_{\mathbf{r}}$, and we shall write $a\xrightarrow[{\mathbf{r}}]{}C_{h}$, iff $b_{h-1}^{\mathbf{r}}\leqslant Ch_{\mathbf{r}}(a)<b_{h}^{\mathbf{r}}$; in particular, in the uncommon case in which $Ch_{\mathbf{r}}(a)=b^{\mathbf{r}}_{p_{\mathbf{r}}}$, then $a$ is assigned to $C_{p_{\mathbf{r}}}$ on $g_{\mathbf{r}}$,
\item $a$ is assigned at most to $C_{h}$ on $g_{\mathbf{r}}$, and we shall write $a\xrightarrow[{\mathbf{r}}]{}C_{\leqslant h}$, iff $Ch_{\mathbf{r}}(a)<b^{\mathbf{r}}_{h}$ (the inequality becomes weak if $h=p_{\mathbf{r}}$),
\item $a$ is assigned at least to $C_{h}$ on $g_{\mathbf{r}}$, and we shall write $a\xrightarrow[{\mathbf{r}}]{}C_{\geqslant h}$, iff $Ch_{\mathbf{r}}(a)\geqslant b^{\mathbf{r}}_{h}$,
\item $a$ is assigned to $[C_{h_1},C_{h_2}]$ on $g_{\mathbf{r}}$, $1\leqslant h_1<h_2\leqslant p_{\mathbf{r}}$, and we shall write $a\xrightarrow[{\mathbf{r}}]{}\left[C_{h_1},C_{h_2}\right]$, iff $b^{\mathbf{r}}_{h_{1}-1}\leqslant Ch_{\mathbf{r}}(a)<b^{\mathbf{r}}_{h_2}$ (the second inequality becomes weak if $h_2=p_{\mathbf{r}}$).
\end{itemize}

To find the parameters necessary to compute the Choquet integral of the alternatives at hand as well as the value of the thresholds delimiting the classes, it can be chosen a direct or an indirect technique \citep{JacquetLagrezeSiskos2001}. In the first case, the DM has to specify all parameters involved in the model (the M\"{o}bius coefficients and the thresholds ${b}^{\mathbf{r}}_{h}$). In the second one, the DM has to provide some preference information in terms of class assignments of some alternative he knows well $\left(a\xrightarrow[{\mathbf{r}}]{}C_{h},\; a\xrightarrow[{\mathbf{r}}]{}C_{\leqslant h},\;a\xrightarrow[{\mathbf{r}}]{}C_{\geqslant h},\; \mbox{or}\;a\xrightarrow[{\mathbf{r}}]{}\left[C_{h_1},C_{h_2}\right]\right)$, pairwise comparison of some alternatives ($a$ is preferred to $b$ on $g_{\mathbf{r}}$ or $a$ and $b$ are indifferent on $g_{\mathbf{r}}$), comparisons between criteria in terms of their importance (``$g_{(\mathbf{r},w_1)}$ is more important than $g_{(\mathbf{r},w_2)}$" or ``$g_{(\mathbf{r},w_1)}$ and $g_{(\mathbf{r},w_2)}$ are equally important", with $g_{(\mathbf{r},w_1)},g_{(\mathbf{r},w_2)}\in G_{\mathbf{r}}^{k}$ for some $\mathbf{r}$ and some $k$) or possible interactions between criteria ($g_{(\mathbf{r},w_1)}$ and $g_{(\mathbf{r},w_2)}$ are positively [negatively] interacting) and parameters compatible with these assignments are therefore inferred in a computational way. After that all preferences provided by the DM are converted into linear inequalities (see below), the set of constraints representing them is the following: 

$$
\left.
\begin{array}{l}
\left.
\begin{array}{l}
\;Ch_{\mathbf{r}}(a)\geqslant b^{\mathbf{r}}_{h-1},\\[2mm]
\;Ch_{\mathbf{r}}(a)-b_{h}^{\mathbf{r}}\leqslant -\varepsilon
\end{array}
\right\}\mbox{if $a\xrightarrow[{\mathbf{r}}]{}C_{h}$}\\[6mm]
\;\;\;Ch_{\mathbf{r}}(a)\geqslant b^{\mathbf{r}}_{h-1},\;\mbox{if $a\xrightarrow[{\mathbf{r}}]{}C_{\geqslant h}$}\\[2mm]
\;\;\;Ch_{\mathbf{r}}(a)-b^{\mathbf{r}}_{h}\leqslant-\varepsilon,\;\mbox{if $a\xrightarrow[{\mathbf{r}}]{}C_{\leqslant h}$}\\[2mm]
\left.
\begin{array}{l}
\;Ch_{\mathbf{r}}(a)\geqslant b^{\mathbf{r}}_{h_1-1},\\[2mm]
\;Ch_{\mathbf{r}}(a)-b_{h_2}^{\mathbf{r}}\leqslant -\varepsilon
\end{array}
\right\}\mbox{if $a\xrightarrow[{\mathbf{r}}]{}[C_{h_1},C_{h_2}]$}\\[6mm]
\left.
\begin{array}{l}
\;\; b^{\mathbf{r}}_{0}=0, \\[2mm]
\;\; b_{p_{\mathbf{r}}}^{\mathbf{r}}=Ch_{\mathbf{r}}\left(opt_{\mathbf{t}_1},\ldots,opt_{\mathbf{t_r}}\right), \\[2mm]
\;\; b^{\mathbf{r}}_{h} \geqslant  b^{\mathbf{r}}_{h-1} + \varepsilon,\; h=1,\ldots,p_{\mathbf{r}}
\end{array}
\right\} \mbox{for all}\;\; g_{\mathbf{r}}\in G\setminus G_E,\\[8mm]
\;\;\;E^{OtherConstrs}\\
\end{array}\\
\right\}{E}_{Ch}^{DM}
$$
\noindent where $\varepsilon$ is an auxiliary variable used to convert the strict inequalities in weak ones\footnote{For example, constraint $x>y$ is converted into $x\leqslant y+\varepsilon.$} and $E^{OtherConstrs}$ is the set of constraints translating the further preference information we specified above and that can be the following:
\begin{itemize}
\item $Ch_{\mathbf{r}}(a)\geqslant Ch_{\mathbf{r}}(b)+\varepsilon$ iff $a$ is preferred to $b$ on $g_{\mathbf{r}}$,
\item $Ch_{\mathbf{r}}(a)=Ch_{\mathbf{r}}(b)$ iff $a$ and $b$ are indifferent on $g_{\mathbf{r}}$\footnote{See \cite{BrankeEtAl2017} for different ways of translating the indifference between alternatives},
\item $\varphi_{\mathbf{r}}\left(\left\{g_{\left(\mathbf{r},w_1\right)}\right\}\right)\geqslant \varphi_{\mathbf{r}}\left(\left\{g_{\left(\mathbf{r},w_2\right)}\right\}\right)+\varepsilon$ iff $g_{\left(\mathbf{r},w_1\right)}$ is more important than $g_{\left(\mathbf{r},w_2\right)}$,
\item $\varphi_{\mathbf{r}}\left(\{g_{\left(\mathbf{r},w_{1}\right)},g_{\left(\mathbf{r},w_{2}\right)}\}\right)\geqslant \varepsilon$ iff $g_{\left(\mathbf{r},w_{1}\right)}$ and $g_{\left(\mathbf{r},w_{2}\right)}$ are positively interacting,
\item $\varphi_{\mathbf{r}}\left(\{g_{\left(\mathbf{r},w_{1}\right)},g_{\left(\mathbf{r},w_{2}\right)}\}\right)\leqslant -\varepsilon$ iff $g_{\left(\mathbf{r},w_{1}\right)}$ and $g_{\left(\mathbf{r},w_{2}\right)}$ are negatively interacting.
\end{itemize}

Denoting by $E^{TC}$ the following set of technical constraints related to the considered preference model

$$
\left.
\begin{array}{l}
\displaystyle\sum_{g_{\mathbf{t}_i}\in G_E}m\left(\{g_{\mathbf{t}_i}\}\right)+\sum_{\{g_{\mathbf{t}_i},g_{\mathbf{t}_j}\}\subseteq G_{E}}m\left(\{g_{\mathbf{t}_i},g_{\mathbf{t}_j}\}\right)=1,\\[2mm]
\forall g_{\mathbf{t}_i}\in G_{E}\;\mbox{and}\;\forall T\subseteq G_{E}\setminus\{g_{\mathbf{t}_i}\}, \;m\left(\{g_{\mathbf{t}_i}\}\right)+\displaystyle\sum_{g_{\mathbf{t}_j}\in T}m\left(\{g_{\mathbf{t}_i},g_{\mathbf{t}_j}\}\right)\geqslant 0
\end{array}
\right\}E_{Ch}^{TC}
$$

\noindent one can check for the existence of at least one instance of the assumed preference model compatible with the preferences provided by the DM solving the following LP problem:

$$
\begin{array}{l}
\varepsilon_{Ch}^{*}=\max\varepsilon, \;\mbox{subject to}, \\[0,5mm]
{E}_{Ch}^{DM}\cup E_{Ch}^{TC}.
\end{array}
$$

\noindent Two cases can occur: 
\begin{itemize}
\item iff ${E}_{Ch}^{DM}\cup E_{Ch}^{TC}$ is feasible and $\varepsilon_{Ch}^{*}>0$, then there exists at least one instance of the considered preference model compatible with the preferences provided by the DM (briefly ``\textit{a compatible model}"),
\item iff ${E}_{Ch}^{DM}\cup E_{Ch}^{TC}$ is infeasible or $\varepsilon_{Ch}^{*}\leqslant 0$, then, there not exists any instance of the assumed preference model compatible with the considered preferences and the cause of this incompatibility can be detected by using one of the procedure proposed in \cite{mousseau2003resolving} and \cite{MousseauDiasFigueira2006}. 
\end{itemize}

\subsection{Looking for a parsimonious representation of the preference information}\label{MinimalSets}
As already underlined in Section \ref{MCHPChoquet}, the use of the 2-additive Choquet integral involves the knowledge of $|E_{G}|+\binom{|E_G|}{2}$ values, one for each elementary criterion and one for each pair of elementary criteria. Anyway, it could be advisable to represent the preferences of the DM by using a preference model as parsimonious as possible in terms of the number of parameters involved. For such a reason, if there exists at least one instance of the 2-additive Choquet integral compatible with the preferences provided by the DM, one can ask which is the minimum number of pairs of interacting criteria necessary to represent these preferences. Therefore, in the following, we shall describe in detail how to compute a minimal set of pairs of interacting criteria such that the 2-additive Choquet integral restores the preferences provided by the DM. In the following, for simplicity, we shall speak of \textit{minimal set of interacting criteria}. Let us observe that the following methodology presents a great interest not only for the introduced method but also for all the MCDA methods in which the underlying preference model is the 2-additive Choquet integral. 

Assuming that the set of constraints $E_{Ch}^{DM}\cup E_{Ch}^{TC}$ is feasible and that $\varepsilon_{Ch}^{*}>0$ (therefore there exists at least one compatible model), and considering a binary variable $\gamma_{\mathbf{t}_i\mathbf{t}_j}\in\{0,1\}$ for each pair of elementary criteria $\{g_{\mathbf{t}_{i}},g_{\mathbf{t}_{j}}\}\subseteq G_{E}$, to check for a minimal set of interacting criteria one has to solve the following MILP problem:

\begin{equation}\label{FirstMinimalSet}
\begin{array}{l}
\;\;\;\;\;\underline{\gamma}_1^{*}=\displaystyle\min\sum_{\{\mathbf{t}_i,\mathbf{t}_j\}\subseteq E_G}\gamma_{\mathbf{t}_i\mathbf{t}_j}, \;\;\mbox{subject to}, \\[7mm]
\left.
\begin{array}{l}
\;\;E_{Ch}^{DM}\cup E_{Ch}^{TC},\\[2mm]
\;\; \delta\leqslant\varepsilon\leqslant\varepsilon_{Ch}^{*},\\[2mm]
\left.
\begin{array}{l}
m\left(\{g_{\mathbf{t}_i},g_{\mathbf{t}_j}\}\right)\geqslant-\gamma_{\mathbf{t}_i\mathbf{t}_j},\\[2mm]
m\left(\{g_{\mathbf{t}_i},g_{\mathbf{t}_j}\}\right)\leqslant\gamma_{\mathbf{t}_i\mathbf{t}_j},\\[2mm]
\end{array}
\right\}\; \mbox{for all}\; \{\mathbf{t}_i,\mathbf{t}_j\}\subseteq E_{G},\\[2mm]
\end{array}
\right\}E_{Ch}^{Min}
\end{array}
\end{equation} 

\noindent where $\delta$ is a very small positive number fixed to ensure that $\varepsilon$ is positive. 

As proved in \cite{grabisch2015exact}, using a 2-additive Choquet integral, the M\"{o}bius parameters of a pair of criteria take values in the interval $\left[-1,1\right]$. Consequently, last two constraints in $E_{Ch}^{Min}$ ensure that if $\gamma_{\mathbf{t}_i\mathbf{t}_j}=0$, then the corresponding M\"{o}bius parameter $m\left(\{g_{\mathbf{t}_i},g_{\mathbf{t}_j}\}\right)$ is null and, therefore, there is not any interaction between $g_{\mathbf{t}_i}$ and $g_{\mathbf{t}_j}$. In the opposite case $\left(\gamma_{\mathbf{t}_i\mathbf{t}_j}=1\right)$, there is a positive or negative interaction between the elementary criteria $g_{\mathbf{t}_i}$ and $g_{\mathbf{t}_j}$. Two cases can be considered: 
\begin{description}
\item[case 1)] $\underline{\gamma}_1^{*}=0$: all pairs of elementary criteria do not present any interaction and, therefore, a weighted sum is able to represent the preferences provided by the DM\footnote{Indeed, the fact that $E^{DM}\cup E^{TC}$ is feasible and $\varepsilon^{*}>0$ implies that there is at least one instance of the 2-additive Choquet integral compatible with the preferences provided by the DM and the weighted sum is a particular case of all Choquet integrals (not only 2-additive ones) since it corresponds to the case in which all interactions are null.},
\item[case 2)] $\underline{\gamma}_1^{*}>0$: a weighted sum is not able to represent the preferences provided by the DM and it is necessary to take into account the interaction between at least one pair of elementary criteria to represent these preferences. In particular, the value $\underline{\gamma}_1^{*}$ is the minimum number of pairs of interacting criteria necessary to represent the provided preferences. 
\end{description}

In \textbf{case 2)}, denoting by $\Gamma_1^{*}$ the vector of binary variables obtained as a solution of the MILP problem (\ref{FirstMinimalSet}), the set $\Gamma_{1}^{Min}=\left\{\gamma_{\mathbf{t}_i\mathbf{t}_j}\in\Gamma_{1}^{*}:\gamma_{\mathbf{t}_i\mathbf{t}_j}=1\right\}$ gives information on the pairs of interacting criteria, while the sign of the corresponding M\"{o}bius parameters $m\left(\{g_{\mathbf{t}_{i}},g_{\mathbf{t}_{j}}\}\right)$ gives information on the type of interaction (positive if $m\left(\{g_{\mathbf{t}_{i}},g_{\mathbf{t}_{j}}\}\right)>0$ or negative if $m\left(\{g_{\mathbf{t}_{i}},g_{\mathbf{t}_{j}}\}\right)<0$). 

In general, more than one set of interacting criteria can represent the provided preferences and, following \cite{MUSA_int}, it could be interesting finding all of them. As explained above, solving the MILP problem (\ref{FirstMinimalSet}), we obtain the first minimal set of interacting criteria $\Gamma^{Min}_{1}$. To check for another possible minimal set of interacting criteria having the same cardinality of $\Gamma^{Min}_{1}$, one has to solve the following MILP problem 

\begin{equation}\label{InteractingProb}
\begin{array}{l}
\;\;\;\underline{\gamma}_2^{*}=\displaystyle\min\sum_{\{\mathbf{t}_i,\mathbf{t}_j\}\subseteq E_{G}}\gamma_{\mathbf{t}_i\mathbf{t}_j}, \;\;\mbox{subject to}, \\[2mm]
\left.
\begin{array}{l}
\displaystyle\sum_{\{\mathbf{t}_i,\mathbf{t}_j\}\subseteq E_{G}}\gamma_{\mathbf{t}_i\mathbf{t}_j}=\underline{\gamma}^{*}_{1},\\[0,6cm]
E_{Ch}^{Min},\\[2mm]
\displaystyle\sum_{\gamma_{\mathbf{t}_i\mathbf{t}_j}\in\Gamma_{1}^{Min}}\gamma_{\mathbf{t}_i\mathbf{t}_j}\leqslant |\Gamma_{1}^{Int}|-1\\[2mm]
\end{array}
\right\}E^{Min_2}_{Ch}
\end{array}
\end{equation}
 
\noindent where the last constraint is included to avoid that the minimal set $\Gamma_{1}^{Min}$ of pairs of interacting criteria is obtained again. If $E^{Min_2}_{Ch}$ is infeasible, then there is not any other minimal set of interacting criteria, while, in the opposite case there exists another minimal set of interacting criteria compatible with the preferences of the DM. Denoting by $\Gamma_{2}^{*}$ the vector of binary variables $\gamma_{\mathbf{t}_i\mathbf{t}_j}$ obtained by solving the MILP problem (\ref{InteractingProb}) and defining $\Gamma_{2}^{Min}=\left\{\gamma_{\mathbf{t}_i\mathbf{t}_j}\in\Gamma_{2}^{*}:\gamma_{\mathbf{t}_i\mathbf{t}_j}=1\right\}$ the subset composed of the binary variables equal to 1, this gives information on the new minimal set of interacting criteria. One can therefore proceed in an iterative way to check for other minimal sets of interacting criteria. Then, if $\Gamma_{k-1}^{Min}$ is the set of binary variables equal to 1 obtained at the iteration $k-1$, the possible $k$-th set of pairs of interacting criteria is obtained by solving the MILP problem 

\begin{equation}\label{InteractingProb_k}
\begin{array}{l}
\;\;\;\underline{\gamma}_k^{*}=\displaystyle\min\sum_{\{\mathbf{t}_i,\mathbf{t}_j\}\subseteq E_{G}}\gamma_{\mathbf{t}_i\mathbf{t}_j}, \;\;\mbox{subject to}, \\[2mm]
\left.
\begin{array}{l}
\displaystyle\sum_{\{\mathbf{t}_i,\mathbf{t}_j\}\subseteq E_{G}}\gamma_{\mathbf{t}_i\mathbf{t}_j}=\underline{\gamma}^{*}_{1},\\[0,6cm]
E_{Ch}^{Min},\\[2mm]
\displaystyle\sum_{\gamma_{\mathbf{t}_i\mathbf{t}_j}\in\Gamma_{1}^{Min}}\gamma_{\mathbf{t}_i\mathbf{t}_j}\leqslant |\Gamma_{1}^{Min}|-1,\\[2mm]
\ldots\;\;\ldots\;\;\ldots \\[2mm]
\displaystyle\sum_{\gamma_{\mathbf{t}_i\mathbf{t}_j}\in\Gamma_{k-1}^{Min}}\gamma_{\mathbf{t}_i\mathbf{t}_j}\leqslant |\Gamma_{k-1}^{Min}|-1.\\[2mm]
\end{array}
\right\}E^{Min_{k}}_{Ch}
\end{array}
\end{equation}
Once all minimal sets of interacting criteria have been obtained, following a basic idea of rough set theory \cite{pawlak1991rough}, one can compute the \textit{core} being the set of pairs of interacting elementary criteria $\{g_{\mathbf{t}_i},g_{\mathbf{t}_j}\}$ such that $\gamma_{\mathbf{t}_i\mathbf{t}_j}\in \Gamma_{1}^{Min}\cap\ldots\cap\Gamma_{k}^{Min}$. The core gives therefore the set of pairs of interacting criteria that are present in all minimal sets and, therefore, are essential to represent the preferences provided by the DM. Let us observe that it is possible that a pair of interacting criteria in the core presents positive interaction in some minimal sets and negative interaction in others as will be shown in the case study described in Section \ref{CaseStudy}. 

\subsection{Non Additive Robust Ordinal Regression applied to the hierarchical and interacting sorting method}\label{NARORsubsec}
If there exist at least one compatible model, in general, there are many of them. All compatible models provide the same recommendations on the reference alternatives in ${A}^{DM}\subseteq{A}$ but they can give different information on the alternatives in $A\setminus A^{DM}$. For such a reason, the choice of only one compatible instance is unjustified or at least arbitrary. Therefore, in this section, we shall show how to apply the Robust Ordinal Regression (ROR; \citealt{GrecoMousseauSlowinski2008}) to the new proposal. ROR provides recommendations on the problem at hand taking into account the whole set of compatible models computing a necessary and a possible preference relation. The extension of ROR to the Choquet integral preference model is known as Non Additive Robust Ordinal Regression (NAROR) and it has been proposed by \cite{angilella2010non}. Considering the hierarchical structure of criteria, the application of NAROR to the proposed method permits to define the following necessary and possible assignments for each alternative $a\in A$, for each non-elementary criterion $g_{\mathbf{r}}$ and for each class $C_h$, $h=1,\ldots,p_{\mathbf{r}}$ (see also \citealt{greco2010multiple}): 

\begin{enumerate}
\item $a$ is necessarily assigned to $C_{h}$ on $g_{\mathbf{r}}$, and we shall write $a\xrightarrow[\mathbf{r}]{N}C_{h}$, iff $b^{\mathbf{r}}_{h-1}\leqslant Ch_{\mathbf{r}}(a)<b^{\mathbf{r}}_{h}$ for all compatible models (the second inequality becomes weak if $h=p_{\mathbf{r}}$),
\item $a$ is possibly assigned to $C_{h}$ on $g_{\mathbf{r}}$, and we shall write $a\xrightarrow[\mathbf{r}]{P}C_{h}$, iff $b^{\mathbf{r}}_{h-1}\leqslant Ch_{\mathbf{r}}(a)<b^{\mathbf{r}}_{h}$ for at least one compatible model (the second inequality becomes weak if $h=p_{\mathbf{r}}$),
\item $a$ is necessarily assigned at least to $C_{h}$ on $g_{\mathbf{r}}$, and we shall write $a\xrightarrow[\mathbf{r}]{N}C_{\geqslant h}$, iff $Ch_{\mathbf{r}}(a)\geqslant b^{\mathbf{r}}_{h-1}$ for all compatible models,
\item $a$ is possibly assigned at least to $C_{h}$ on $g_{\mathbf{r}}$, and we shall write $a\xrightarrow[\mathbf{r}]{P}C_{\geqslant h}$, iff $Ch_{\mathbf{r}}(a)\geqslant b^{\mathbf{r}}_{h-1}$ for at least one compatible model, 
\item $a$ is necessarily assigned at most to $C_{h}$ on $g_{\mathbf{r}}$, and we shall write $a\xrightarrow[\mathbf{r}]{N}C_{\leqslant h}$, iff $Ch_{\mathbf{r}}(a)<b^{\mathbf{r}}_{h}$ for all compatible models (the inequality becomes weak if $h=p_{\mathbf{r}}$),
\item $a$ is possibly assigned at most to $C_{h}$ on $g_{\mathbf{r}}$, and we shall write $a\xrightarrow[\mathbf{r}]{P}C_{\leqslant h}$, iff $Ch_{\mathbf{r}}(a)<b^{\mathbf{r}}_{h}$ for at least one compatible model (the inequality becomes weak if $h=p_{\mathbf{r}}$).
\end{enumerate}

From the definitions of necessary and possible assignments above, it is clear that for each $a\in A$ and for each non-elementary criterion $g_{\mathbf{r}}$, the following two properties must hold\footnote{Analogous properties for a hierarchical sorting in case the assumed preference model is an additive value function can be found in \cite{CorrenteEtAl2017a}.}:
\begin{itemize}
	\item either $a$ is necessarily assigned at least to $C_h$ on $g_{\mathbf{r}}$ $\left(a\xrightarrow[\mathbf{r}]{N}C_{\geqslant h}\right)$  or $a$ is possibly assigned at most to $C_{h-1}$ on $g_{\mathbf{r}}$ $\left(a\xrightarrow[\mathbf{r}]{P}C_{\leqslant h-1}\right)$, with $h\in\left\{2,\ldots,p_{\mathbf{r}}\right\}$,
	\item either $a$ is necessarily assigned at most to $C_h$ on $g_{\mathbf{r}}$ $\left(a\xrightarrow[\mathbf{r}]{N}C_{\leqslant h}\right)$  or $a$ is possibly assigned at least to $C_{h+1}$ on $g_{\mathbf{r}}$ $\left(a\xrightarrow[\mathbf{r}]{P}C_{\geqslant h+1}\right)$, with $h\in\left\{1,\ldots,p_{\mathbf{r}}-1\right\}$.
\end{itemize}

From a computational point of view, for each $a\in A$, for each non-elementary criterion $g_{\mathbf{r}}$ and for each $C_h$, $h=1,\ldots,p_{\mathbf{r}}$, the necessary and possible assignments above can be obtained as follows: 

\begin{enumerate}\addtocounter{enumi}{1}
\item Defined $E^{P}_{h}=E^{DM}_{Ch} \cup E^{TC}_{Ch}\cup\{Ch_{\mathbf{r}}(a)\geqslant b_{h-1}^{\mathbf{r}};\;Ch_{\mathbf{r}}(a)+\varepsilon\leqslant b_{h}^{\mathbf{r}}\}$, $a\xrightarrow[\mathbf{r}]{P}C_{\leqslant h}$ iff $E^{P}_{h}$ is feasible and $\varepsilon_{h}^{P}>0$, where $\varepsilon_{h}^{P}=\max\varepsilon$ subject to $E^{P}_{h}$\footnote{$E^{P}_{h}=E^{DM}_{Ch} \cup E^{TC}_{Ch}\cup\{Ch_{\mathbf{r}}(a)\geqslant b_{h-1}^{\mathbf{r}};\;Ch_{\mathbf{r}}(a)\leqslant b_{h}^{\mathbf{r}}\}$ if $h=p_{\mathbf{r}}$.},
\item Defined $E^{N}_{\geqslant h}=E_{Ch}^{DM}\cup E_{Ch}^{TC}\cup\left\{Ch_{\mathbf{r}}(a)\leqslant b_{h-1}^{\mathbf{r}}-\varepsilon\right\}$, $a\xrightarrow[\mathbf{r}]{N}C_{\geqslant h}$ iff $E^{N}_{\geqslant h}$ is infeasible or $\varepsilon_{\geqslant h}^{N}\leqslant 0$, where $\varepsilon_{\geqslant h}^{N}=\max\varepsilon$ subject to $E^{N}_{\geqslant h}$,
\item Defined $E^{P}_{\geqslant h}=E_{Ch}^{DM}\cup E_{Ch}^{TC}\cup\left\{Ch_{\mathbf{r}}(a)\geqslant b_{h-1}^{\mathbf{r}}\right\}$, $a\xrightarrow[\mathbf{r}]{P}C_{\geqslant h}$ iff $E^{P}_{\geqslant h}$ is feasible and $\varepsilon_{\geqslant h}^{P}>0$, where $\varepsilon_{\geqslant h}^{P}=\max\varepsilon$ subject to $E^{P}_{\geqslant h}$,
\item Defined $E^{N}_{\leqslant h}=E_{Ch}^{DM}\cup E_{Ch}^{TC}\cup\left\{Ch_{\mathbf{r}}(a)\geqslant b_{h}^{\mathbf{r}}\right\}$, $a\xrightarrow[\mathbf{r}]{N}C_{\leqslant h}$ iff $E^{N}_{\leqslant h}$ is infeasible or $\varepsilon_{\leqslant h}^{N}\leqslant 0$, where $\varepsilon_{\leqslant h}^{N}=\max\varepsilon$ subject to $E^{N}_{\leqslant h}$\footnote{$E^{N}_{\leqslant h}=E_{Ch}^{DM}\cup E_{Ch}^{TC}\cup\left\{Ch_{\mathbf{r}}(a)\geqslant b_{h}^{\mathbf{r}}+\varepsilon\right\}$ if $h=p_{\mathbf{r}}$.},
\item Defined $E^{P}_{\leqslant h}=E_{Ch}^{DM}\cup E_{Ch}^{TC}\cup\left\{Ch_{\mathbf{r}}(a)\leqslant b_{h}^{\mathbf{r}}-\varepsilon\right\}$, $a\xrightarrow[\mathbf{r}]{P}C_{\leqslant h}$ iff $E^{P}_{\leqslant h}$ is feasible and $\varepsilon_{\leqslant h}^{P}> 0$, where $\varepsilon_{\leqslant h}^{P}=\max\varepsilon$ subject to $E^{P}_{\leqslant h}$\footnote{$E^{P}_{\leqslant h}=E_{Ch}^{DM}\cup E_{Ch}^{TC}\cup\left\{Ch_{\mathbf{r}}(a)\leqslant b_{h}^{\mathbf{r}}\right\}$ if $h=p_{\mathbf{r}}$.}.
\end{enumerate}

To avoid to solve a MILP problem, the necessary assignment of an alternative $a$ to a single class $C_h$ on criterion $g_{\mathbf{r}}$ can be computed on the basis of $a\xrightarrow[\mathbf{r}]{N}C_{\geqslant h}$ and $a\xrightarrow[\mathbf{r}]{N}C_{\leqslant h}$. In particular, $a\xrightarrow[\mathbf{r}]{N}C_{h}$ iff $a\xrightarrow[\mathbf{r}]{N}C_{\geqslant h}$ and $a\xrightarrow[\mathbf{r}]{N}C_{\leqslant h}$.

\subsection{Stochastic Multicriteria Acceptability Analysis applied to the hierarchical and interacting sorting method}\label{SMAAsubsec}
Looking at the assignments provided by the application of the NAROR, it is very common that an alternative is possibly assigned to more than one class with respect to a particular criterion in the hierarchy as well as at comprehensive level. For such a reason, it is useful to have an estimate of the number of compatible instances, with respect to the whole set of compatible instances, giving a certain assignment. To this aim, we shall propose the application of the SMAA methodology that has been introduced at first in \cite{LahdelmaHokkanenSalminen1998} to deal with ranking problems in which the underlying preference model was a weighted sum (see \citealt{PelissariEtAl2019} for a recent survey on SMAA). After that, SMAA has been further extended and applied to sorting problems \citep{milos_tommi_2} as well as to the Choquet integral preference model \citep{angilella2015,angilella2016MCHP}. \\
The application of SMAA to sorting problems permits to explore the whole set of models compatible with the preferences provided by the DM giving recommendations in terms of frequency of assignment of an alternative to a certain class with respect to a particular criterion in the hierarchy. 

Denoting by $\left(\mbox{\textbf{\textit{Ch}}},{\cal B}\right)$ the whole set of compatible models (that is vectors $\left(Ch,\mathbf{b}\right)$ such that the M\"{o}bius transform of a capacity in $Ch$ and the threshold values in $\mathbf{b}$ are compatible with the preferences provided by the DM), for each $a\in A$, for each non-elementary criterion $g_{\mathbf{r}}$ and for each class $C_h$, $h=1,\ldots,p_{\mathbf{r}}$, the class acceptability index $CAI(a,h,\mathbf{r})$ is computed as follows \citep{milos_tommi_2}:

$$
CAI(a,h,\mathbf{r})=\int_{\left(Ch,\mathbf{b}\right)\in\left(\mbox{\textbf{\textit{Ch}}},{\cal B}\right)} CMF(Ch,\mathbf{b},a,h,\mathbf{r}) \;dP\left(Ch,\mathbf{b}\right).
$$

\noindent where $P\left(Ch,\mathbf{b}\right)$ is a probability measure on $\left(\mbox{\textbf{\textit{Ch}}},{\cal B}\right)$ which is usually assumed as uniform or inferred from some DM preference information (see \cite{}). $CAI(a,h,\mathbf{r})\in\left[0,1\right]$ is the frequency with which $a$ is assigned to $C_{h}$ on $g_{\mathbf{r}}$, and $CMF(Ch,\mathbf{b},a,h,\mathbf{r})$ is a class membership function being equal to $1$ iff $Ch_{\mathbf{r}}(a)\in[b^{\mathbf{r}}_{h-1},b^{\mathbf{r}}_{h}[$ and $0$ otherwise. Of course, the same membership function is equal to $1$ if $Ch_{\mathbf{r}}(a)\in[b^{\mathbf{r}}_{h-1},b^{\mathbf{r}}_{h}]$ in case $h=p_{\mathbf{r}}$.

From a technical point of view, the integral above can be computed by the application of Monte Carlo methods (see \citealt{CorrenteGrecoSlowinski2019} for a more detailed description of the application of Monte Carlo methods to compute multidimensional integrals in case the Choquet integral is assumed as preference model).

To further support the DM in the choice of a final classification, let us show how to integrate the class acceptability indices computed above. Even if some approaches to aggregate the results of the application of the SMAA for ranking problems have been presented (for example, \citealt{kadzinski2016scoring} and \citealt{vetschera2017deriving}), to the best of our knowledge a similar approach for sorting problems has not been proposed. For such a reason, in the following lines we shall briefly describe a procedure which aims to minimize misleading classifications considering the behavior of the DM. In this perspective, for each non-elementary criterion $g_{\mathbf{r}}$ we consider a \textit{loss function} \citep{gneiting2007strictly,savage1971elicitation,schervish1989general} defined as follow
\begin{equation}\label{LossFunc}
{L}(\mathbf{y})=\displaystyle\sum_{a \in A} \sum_{h=1}^{p_{\mathbf{r}}}y_{a,h}\sum_{k\neq h}d(C_{k},C_{h})CAI(a,k,\mathbf{r}),
\end{equation}
\noindent where
\begin{itemize}
\item $\mathbf{y}=[{y}_{a,h},\; a\in A,\;h=1,\ldots,{p}_{\mathbf{r}}],$ with ${y}_{a,h}=1$ if alternative $a$ is assigned to the class ${C}_{h}$ and ${y}_{a,h}=0$ otherwise;
\item $d({C}_{k},{C}_{h})$ is a function which assigns a weight to $CAI(a,k,\mathbf{r})$ depending on the distance between the classes $C_k$ and $C_h$. As in \cite{rademaker2009loss}, $d(C_{k},C_{h})$ has not to be a metric but it has to satisfy (i) $d({C}_{k},{C}_{h})\geqslant 0$ for all $h,k\in\{1,\ldots,p_{\mathbf{r}}\}$; (ii) $d({C}_{k},{C}_{h})=0\Leftrightarrow{C}_{k}={C}_{h};$ (iii) $d({C}_{i},{C}_{j})\leqslant d({C}_{i},{C}_{k})$ and $d({C}_{j},{C}_{k})\leqslant d({C}_{i},{C}_{k})$ for all ${C}_{i}<{C}_{j}<{C}_{k}.$
\end{itemize}

Clearly, since the choice of the distance $d(C_{k},C_{h})$ in eq. (\ref{LossFunc}) affects the assignments, this procedure requires to involve the DM into the process. For example, one DM could decide that the loss in assigning an alternative to a wrong class is marginally decreasing and, therefore, $d({C}_{h},{C}_{k})=\sqrt{|h-k|}$ or one can consider that the same loss is marginally increasing and therefore, setting a distance as $d({C}_{h},{C}_{k})=({h-k})^{2}$.

Solving the following MILP problem, 

\begin{equation}\label{LossFunctionProblem}
\begin{array}{l}
\displaystyle L^{*}=\min_{\mathbf{y}}L(\mathbf{y}), \;\;\mbox{subject to}\\
\displaystyle\sum_{h=1}^{p_{\mathbf{r}}}y_{a,h}=1, \;\;\mbox{for all}\;a\in A
\end{array}
\end{equation}

\noindent we get a single assignment summarizing the class acceptability indices previously computed. One can wonder if there exists another assignment presenting the same $L^{*}$ value. Therefore, denoting by $\mathbf{y}^{*}$ the vector of binary variables obtained as solution of the previous MILP problem, another assignment can be obtained by solving $(\ref{LossFunctionProblem})$ with the addition of the following constraint

$$
\displaystyle\sum_{y^{*}_{a,h}\in y^{*}:\; y^{*}_{a,h}=1}y^{*}_{a,h}\leqslant L^{*}-1, 
$$

\noindent avoiding to find, again, the same solution obtained solving the previous MILP problem. If the new MILP problem is feasible, then there is another possible assignment and others can be found by proceeding in an iterative way as shown above. In the opposite case, the found assignment is the unique one. 

\section{Comparison with other methods dealing with sorting problems in presence of interacting criteria}\label{ComparisonSection}
As already stated in the introduction, very few contributions dealt with the use of the Choquet integral preference model for sorting problems and all of them differ with respect to our proposal in several aspects. In \cite{Roubens2001}, the authors propose the application of the Choquet integral preference model for sorting problems. The method is articulated in two parts. In the first part, for each alternative and for each criterion a net score is computed. It is given by the difference between the number of alternatives to which the alternative is preferred on the criterion and the number of alternatives preferred to it on the same criterion. This pre-scoring phase is performed to put all evaluations on the same scale since this is essential for the application of the Choquet integral preference model. However, we are not dealing with this problem in the paper. Moreover, other methods can also be used to put all evaluations on the same scale such as the one proposed in \cite{angilella2015} or the one proposed in \cite{GrecoEtAl2018} and that will be used in the case study in the next section. The second phase is related, instead, to the use of the Choquet integral to represent a complete classification provided by the DM or a partial classification regarding some alternatives he knows well. In particular, in the second case (being the most similar to our approach), the DM provides some precise classifications that are translated in linear constraints so that if alternatives $a$ and $b$ are assigned to classes $C_{h_1}$ and $C_{h_2}$, respectively, with $h_1>h_{2},$ then $Ch(a)$ has to be greater than $Ch(b)$. Then, by solving an LP problem one checks for the existence of an instance of the Choquet integral compatible with these preferences. The thresholds delimiting the classes are not variables of the problem but they are defined ex-post so that $b_{h-1}=\min\left\{Ch(a): a\in A \;\mbox{and}\; a\xrightarrow[\mbox{DM}]{}C_h\right\}$ that is, $b_{h-1}$ is equal to the minimum value of the Choquet integral of all alternatives assigned by the DM to the class $C_h$. The same thresholds are then used to assign all the other alternatives which are not provided as reference examples by the DM. 

As mentioned above, our proposal is different from \cite{Roubens2001} for several reasons. First of all, in our approach the thresholds are considered as variables and are not defined ex-post from the assignments provided by the DM. Indeed, in \cite{Roubens2001}, the assignment of non-reference alternatives is strictly dependent on which alternatives are provided as reference examples from the DM. Let us explain this point more in detail. Let us assume that the DM is dealing with a simple problem in which 6 alternatives $\{a_{1},\ldots,a_{6}\}$ have to be assigned to 2 different classes $\{C_{1},C_{2}\}$ and that, with respect to his preference system, $a_{1}\succ a_{2}\succ \ldots\succ a_{6}$ (where ${a}_{i}\succ {a}_{k}$ means that ${a}_{i}$ is preferred to ${a}_{k}$) with the first three alternatives $(a_{1},a_{2},a_{3})$ classified in $C_{2}$ and the last three $(a_{4},a_{5},a_{6})$ classified in $C_{1}$ since he retains the first three alternatives as good and the last three as bad. Let also suppose that there exist at least one capacity such that the Choquet integral is able to restore these preferences. If the DM provides information regarding alternatives ${a}_{1}$ and ${a}_{4}$ only, the method proposed in \cite{Roubens2001} would assign ${a}_{1}$ to ${C}_{2}$ putting all the others in class ${C}_{1}$. This because ${b}_{1}=Ch({a}_{1})$ and all the other alternatives have a Choquet evaluation lower than ${b}_{1}$. In such case, the hypothetical classification of the DM would not be restored. Instead, in our approach the thresholds are not fixed. We consider all the thresholds compatible with the preference information provided by the DM. The assignment of only ${a}_{1}$ to ${C}_{2}$ does not avoid that other alternatives could be assigned to ${C}_{2}$ as well. 


The other differences are mainly due to the fact that we are taking into account a hierarchical structure of criteria permitting, therefore, to get alternatives' assignments not only at comprehensive level but also considering a particular node of the hierarchy and, moreover, the fact that we are considering robustness concerns by using the NAROR and the SMAA methodologies. Indeed, while in \cite{Roubens2001} the proposed procedure gives a single capacity for which the Choquet integral is compatible with the preferences provided by the DM and the same capacity is therefore used to perform the assignment of the remaining alternatives, in our approach we are taking into account all the capacities compatible with the provided preferences giving, consequently, more stable results. \\
Let us observe that other contributions developing the procedure introduced in \cite{Roubens2001} have been proposed in \cite{MarichalMeyerRoubens2005,MarichalRoubens2001,MeyerRoubens2005}. In particular, \cite{MeyerRoubens2005} introduced the TOMASO software implementing the procedure previously described, while \cite{MarichalMeyerRoubens2005} studied the potentialities of the same method by means of two examples.  

Finally, another contribution dealing with the use of the Choquet integral preference model (not only the 2-additive one) for sorting problems, is due to \cite{BenabbouPernyViappiani2017}. Anyway, in that paper, the authors deal not only with sorting but also with ranking and choice problems from an artificial intelligence point of view. They propose several procedures to reduce, iteration by iteration, the space of capacities compatible with the preferences provided by the DM\footnote{The interested reader can also have a look at \cite{BrankeEtAl2017} and \cite{CiomekKadzinskiTervonen2016}.}. 

\section{Case study}\label{CaseStudy}
\begin{figure}[!htb]
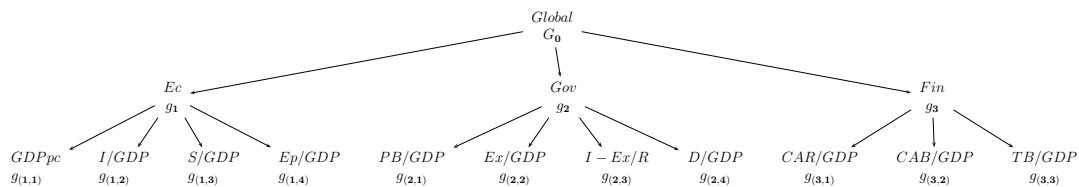

\centering
\caption{Hierarchical structure of criteria considered in the case study\label{Hierarchy}}
\resizebox{0.8\textwidth}{!}{
\pstree[nodesep=2pt]{\TR{$\begin{array}{c}Global\\ G_{\mathbf{0}}\end{array}$}}{
         \pstree{\TR{$\begin{array}{c}Ec\\ g_{\mathbf{1}}\end{array}$}}{
        \TR{$\begin{array}{l}GDPpc\\ g_{\mathbf{(1,1)}}\end{array}$}				
				\TR{$\begin{array}{l}I/GDP\\ g_{\mathbf{(1,2)}}\end{array}$}
				\TR{$\begin{array}{l}S/GDP\\ g_{\mathbf{(1,3)}}\end{array}$}
				\TR{$\begin{array}{l}Ep/GDP\\ g_{\mathbf{(1,4)}}\end{array}$}
        }
        \pstree{ \TR{$\begin{array}{c}Gov\\ g_{\mathbf{2}}\end{array}$}}{
				\TR{$\begin{array}{c}PB/GDP\\ g_{\mathbf{(2,1)}}\end{array}$}				
				\TR{$\begin{array}{c}Ex/GDP\\ g_{\mathbf{(2,2)}}\end{array}$}
				\TR{$\begin{array}{c}I-Ex/R\\ g_{\mathbf{(2,3)}}\end{array}$}
				\TR{$\begin{array}{c}D/GDP\\ g_{\mathbf{(2,4)}}\end{array}$}
				}
				\pstree{ \TR{$\begin{array}{c}Fin\\ g_{\mathbf{3}}\end{array}$}}{
				\TR{$\begin{array}{c}CAR/GDP\\ g_{\mathbf{(3,1)}}\end{array}$}				
				\TR{$\begin{array}{c}CAB/GDP\\ g_{\mathbf{(3,2)}}\end{array}$}
				\TR{$\begin{array}{c}TB/GDP\\ g_{\mathbf{(3,3)}}\end{array}$}
				}
				}}
\end{figure}

In this section, we shall apply the proposed methodology to sovereign ratings of the countries belonging to the European Union (EU)\footnote{\url{https://europa.eu/european-union/index_en}. For the time being, the United Kingdom remains a full member of the EU and rights and obligations continue to fully apply in and to the UK.}. The 28 Countries are evaluated with respect to three macro-criteria and 11 elementary criteria descending from them described in Table \ref{glossary} and articulated as shown in Fig. \ref{Hierarchy}:
\begin{itemize}
\item Economic (Ec) $(g_{\mathbf{1}})$:
\begin{itemize}
\item GDP per capita (GDPc) $(g_{\mathbf{(1,1)}})$, 
\item Investment/GDP (I/GDP) $(g_{\mathbf{(1,2)}})$, 
\item Savings/GDP (S/GDP) $(g_{\mathbf{(1,3)}})$,  
\item Exports/GDP (Ep/GDP) $(g_{\mathbf{(1,4)}})$, 
\end{itemize} 
\item Governmental (Gov) $(g_{\mathbf{2}})$:
\begin{itemize}
\item Primary balance/GDP (PB/GDP) $(g_{\mathbf{(2,1)}})$, 
\item Expenditure/GDP (Ex/GDP) $(g_{\mathbf{(2,2)}})$,  
\item Interest/revenues (I-Ex/R) $(g_{\mathbf{(2,3)}})$,  
\item Net debt/GDP (D/GDP) $(g_{\mathbf{(2,4)}})$, 
\end{itemize}
\item Financial (Fin) $(g_{\mathbf{3}})$: 
\begin{itemize}
\item Current account receipts/GDP (CAR/GDP) $(g_{\mathbf{(3,1)}})$,  
\item Current account balance/GDP (CAB/GDP) $(g_{\mathbf{(3,2)}})$, 
\item Trade balance/GDP (TB/GDP) $(g_{\mathbf{(3,3)}})$. 
\end{itemize}
\end{itemize}

\begin{table}[!htb]
\centering
\caption{Description of the considered criteria\label{glossary}}
\resizebox{1\textwidth}{!}{
		\begin{tabular}{lll}
		\hline
		\hline
		& \textbf{Elementary criterion} & \textbf{Description} \\
		\hline
		\hline
		$g_{\mathbf{1}}$ & {Economic Data} (Ec) & Economic assessments key indicators \\[1mm]
		$g_{\mathbf{2}}$ & {Government Data} (Gov) & Government assessment key indicators \\[1mm]
		$g_{\mathbf{3}}$ & {Balance of Payment Data} (Fin) & External assessment key indicators \\
		
		\hline 
		\hline		
		$g_{(\mathbf{1,1})}$ & {GDP per capita (GDPc)} & Total US dollar market value of goods and services produced by resident factors of production, \\[1mm]
		& & divided by population \\[1mm]
		
		$g_{(\mathbf{1,2})}$ & {Investment/GDP (I/GDP)} & Expenditure on capital goods including plant, equipment, and housing, plus the change in \\[1mm]
		& & inventories, as a percent of GDP \\[1mm]
		
    $g_{(\mathbf{1,3})}$ & {Savings/GDP (S/GDP)} & Investment plus the current account surplus (deficit), as a percent of GDP \\[1mm]
		
		$g_{(\mathbf{1,4})}$ & {Exports/GDP (Ep/GDP)} & Exports of goods and services as a share of GDP \\[1mm]
		
		\hline
		\hline
		
		$g_{(\mathbf{2,1})}$ & {Primary balance/GDP (PB/ GDP)} & Surplus (deficit) plus interest payments on general government debt, as a percent of GDP \\[1mm]
		
		$g_{(\mathbf{2,2})}$ & {Expenditure/GDP (Ex/GDP)} & General government recurrent expenditure, for purposes such as salaries, goods for immediate \\[1mm]
		& & consumption, interest and other transfers, plus capital spending that increases the value of  \\[1mm]
		& & general government physical assets, plus net lending (where a subsidy is involved), as a percent of GDP\\[1mm]
		
		$g_{(\mathbf{2,3})}$ & {Interest/revenues (I-Ex/R)} & Interest payments on general government debt, as a percent of general government revenues \\[1mm]
		
		$g_{(\mathbf{2,4})}$ & {Net debt/GDP (D/GDP)} & Gross debt minus general government financial assets \\[1mm]
		
		\hline
		\hline
		
		$g_{(\mathbf{3,1})}$ & {Current account receipts/GDP (CAR/GDP)} & CAR include proceeds from exports of goods and services, factor income earned by \\[1mm]
		& & residents from nonresidents, and official and private transfers to residents from nonresidents,  \\[1mm]
		& & as a percent of GDP \\[1mm]
		
		$g_{(\mathbf{3,2})}$ & {Current account balance/GDP (CAB/GDP)} & Exports of goods and services minus imports of the same plus net factor income \\[1mm]
		& & plus official and private net transfers, as a percent of GDP \\[1mm]
		
		$g_{(\mathbf{3,3})}$ & {Trade balance/GDP (TB/ GDP)} & Receipts from exports of goods minus payments for imports of goods as a percentage of GDP \\[1mm]
		\hline
		\hline
	  \end{tabular}
		}
\end{table}

The evaluations of the 28 countries on the 11 elementary criteria are shown in Table \ref{EvaluationsTable}. Data have been collected from Standard \& Poor's Global Market Intelligence which published them on December 14th 2017\footnote{\url{https://www.capitaliq.com/CIQDotNet/CreditResearch/RenderArticle.aspx?articleId=1969912\&SctArtId=446671\&from=CM\&nsl_code=LIME\&sourceObjectId=10373140\&sourceRevId=1\&fee_ind=N\&exp_date=20280116-20:50:00}}.

\begin{table}[htbp]
  \centering
  \caption{Performances of the 28 Countries on the 11 elementary criteria taken into account and their rating}
	\resizebox{1\textwidth}{!}{
    \begin{tabular}{l|cccc|cccc|ccc|c|c}
		\hline
		\hline
          & \multicolumn{1}{c}{$g_{\mathbf{(1,1)}}$} & \multicolumn{1}{c}{$g_{\mathbf{(1,2)}}$} & \multicolumn{1}{c}{$g_{\mathbf{(1,3)}}$} & $g_{\mathbf{(1,4)}}$ & \multicolumn{1}{c}{$g_{\mathbf{(2,1)}}$} & \multicolumn{1}{c}{$g_{\mathbf{(2,2)}}$} & \multicolumn{1}{c}{$g_{\mathbf{(2,3)}}$} & $g_{\mathbf{(2,4)}}$ & $g_{\mathbf{(3,1)}}$ & $g_{\mathbf{(3,2)}}$ & $g_{\mathbf{(3,3)}}$ & S\&P rating & Classification\\
		\hline			
		\hline
    Austria & 47801,33 & 24.3  & 26.7  & 54.1  & 0.86  & 50.60 & 3.95  & 73.31 & 62.30 & 2.38  & 0.45 & AA+ & $C_4$ \\
    Belgium & 43762,16 & 24.1  & 24.1  & 85.7  & 1.11  & 52.50 & 5.12  & 97.35 & 99.17 & -0.02 & 0.20 & AA & $C_4$ \\
    Bulgaria & 8058,38 & 19.3  & 22.5  & 64.7  & 0.85  & 34.92 & 2.44  & 15.62 & 71.61 & 3.25  & -3.24 & BBB- & $C_2$ \\
    Cyprus & 25140,98 & 19.3  & 11.9  & 64.8  & 2.89  & 38.00 & 6.20  & 91.17 & 90.47 & -7.37 & -21.76 & BB+ & $C_1$ \\
    Croatia & 13122,86 & 20.2  & 23.2  & 51.3  & 1.93  & 46.10 & 6.73  & 74.68 & 58.65 & 3.07  & -17.11 & BB & $C_1$\\
    Denmark & 56272,84 & 20.8  & 28.3  & 54.2  & 0.08  & 53.20 & 2.45  & 23.95 & 64.86 & 7.47  & 5.64 & AAA & $C_4$ \\
    Estonia & 19453,92 & 24.5  & 25.5  & 77.2  & 0.15  & 40.90 & 0.12  & -3.28 & 84.55 & 0.97  & -4.42 & AA- & $C_4$ \\
    Finland & 45860,33 & 22.6  & 21.3  & 36.1  & -0.20 & 55.30 & 2.04  & 23.82 & 43.53 & -1.23 & 0.40 & AA+ & $C_4$ \\
    France & 38604,63 & 22.8  & 21.6  & 30.0  & -1.06 & 56.10 & 3.46  & 90.45 & 38.71 & -1.14 & -1.27 & AA & $C_4$ \\
    Germany & 45182,19 & 20.2  & 27.7  & 46.1  & 1.47  & 44.70 & 2.82  & 59.84 & 53.48 & 7.53  & 8.10 & AAA & $C_4$ \\
    Greece & 18912,07 & 10.6  & 10.1  & 31.3  & 1.80  & 50.80 & 5.61  & 167.65 & 32.93 & -0.47 & -9.33 & B- & $C_1$ \\
    Ireland & 68030,49 & 32.8  & 33.0  & 128.7 & 1.55  & 27.90 & 7.88  & 61.63 & 153.47 & 0.15  & 41.37 & A+ & $C_3$ \\
    Italy & 32077,55 & 17.2  & 19.5  & 31.2  & 1.70  & 49.30 & 8.06  & 122.74 & 35.40 & 2.35  & 3.16 & BBB & $C_2$ \\
    Latvia & 15570,01 & 21.6  & 21.0  & 59.5  & 0.55  & 37.70 & 3.37  & 31.24 & 68.80 & -0.63 & -9.53 & A- & $C_3$ \\
    Lithuania & 16197,33 & 17.7  & 17.0  & 76.1  & 0.52  & 35.70 & 3.48  & 34.88 & 81.23 & -0.77 & -5.51 & A- & $C_3$ \\
    Luxembourg & 106866,09 & 18.4  & 23.0  & 223.2 & 0.89  & 42.90 & 0.66  & -15.46 & 557.37 & 4.53  & -6.61 & AAA & $C_4$ \\
    Malta & 27295,35 & 22.2  & 31.7  & 136.0 & 2.71  & 37.80 & 4.94  & 44.48 & 236.48 & 9.48  & -16.95 & A- &  $C_3$ \\
    Netherlands & 48474,25 & 20.6  & 29.4  & 85.1  & 1.72  & 42.90 & 2.12  & 53.50 & 119.27 & 8.79  & 12.36 & AAA & $C_4$ \\
    Poland & 13248,49 & 20.3  & 19.5  & 53.8  & 0.27  & 41.10 & 4.99  & 49.46 & 57.37 & -0.79 & -0.30 & BBB+ & $C_2$ \\
    Portugal & 21118,94 & 16.2  & 16.8  & 42.9  & 2.41  & 45.00 & 8.73  & 116.80 & 50.83 & 0.64  & -5.41 & BBB- & $C_2$ \\
    U.K.  & 39493,34 & 17.6  & 12.3  & 29.1  & 0.07  & 40.50 & 6.76  & 83.43 & 36.71 & -5.32 & -6.84 & AA & $C_4$ \\
    Czech Rep. & 20161,53 & 26.5  & 28.3  & 76.1  & 1.63  & 39.20 & 2.07  & 29.40 & 80.70 & 1.80  & 5.52 & AA- & $C_4$ \\
    Romania & 10289,76 & 24.5  & 20.9  & 41.1  & -1.68 & 37.10 & 4.19  & 30.92 & 45.85 & -3.55 & -6.79 & BBB- & $C_2$ \\
    Slovakia & 17532,04 & 23.1  & 22.4  & 94.6  & -0.23 & 41.90 & 4.16  & 46.32 & 96.90 & -0.67 & 2.01 & A+ & $C_3$ \\
    Slovenia & 23504,60 & 19.1  & 24.3  & 79.1  & 1.54  & 44.50 & 5.59  & 58.77 & 84.06 & 5.20  & 3.46 & A+ & $C_3$ \\
    Spain & 28288,91 & 20.9  & 22.6  & 34.1  & -0.59 & 41.10 & 6.61  & 88.92 & 40.09 & 1.69  & -2.06 & BBB+ & $C_2$ \\
    Sweden & 54215,57 & 25.5  & 30.0  & 44.6  & 1.22  & 49.30 & 1.05  & 25.08 & 54.01 & 4.46  & 2.03 & AAA & $C_4$ \\
    Hungary & 13372,27 & 20.8  & 24.8  & 88.6  & 0.46  & 51.50 & 6.03  & 68.62 & 94.13 & 4.01  & 3.08 & BBB- & $C_2$ \\
		\hline
		\hline
    Pref. Direction & \multicolumn{1}{c}{Inc.} & \multicolumn{1}{c}{Inc.} & \multicolumn{1}{c}{Inc.} & Inc. & \multicolumn{1}{c}{Inc.} & \multicolumn{1}{c}{Inc.} & \multicolumn{1}{c}{Dec.} & Dec. & \multicolumn{1}{c}{Inc.} & \multicolumn{1}{c}{Inc.} & \multicolumn{1}{c}{Inc.} \\
		\hline
		\hline
    \end{tabular}%
		}
  \label{EvaluationsTable}%
\end{table}%

Since, as explained in Section \ref{MCHPChoquet}, the application of the Choquet integral implies that all evaluations are expressed on the same scale, we apply the following normalization technique used in \cite{GrecoEtAl2018}. For each elementary criterion $g_{\mathbf{t}}$, let us compute the mean $(M_{\mathbf{t}})$ and the standard deviation $(s_{\mathbf{t}})$ of the evaluations in Table \ref{EvaluationsTable} as follows: 

$$
M_{\mathbf{t}}=\displaystyle\frac{1}{|A|}\sum_{a\in A}g_{\mathbf{t}}(a)\;\;\;\;\;\mbox{and}\;\;\;\;\; s_{\mathbf{t}}=\sqrt{\frac{\displaystyle\sum_{a\in A}\left(g_{\mathbf{t}}(a)-M_{\mathbf{t}}\right)^2}{|A|}}.
$$

\noindent Then, for each $a\in A$ and for each $\mathbf{t}\in E_G$, the $z$-score $g_{\mathbf{t}}^{z}(a)$ is obtained: 

$$
g_{\mathbf{t}}^{z}(a)=\frac{g_{\mathbf{t}}(a)-M_{\mathbf{t}}}{s_{\mathbf{t}}}. 
$$

\noindent Finally, the normalized evaluation $\overline{g}_{\mathbf{t}}(a)$ of $a\in A$ on elementary criterion $g_{\mathbf{t}}$ is computed as 

$$
\overline{g}_{\mathbf{t}}(a)=
\left\{
\begin{array}{lll}
0 & \mbox{if} & g_{\mathbf{t}}(a)\leqslant M_{\mathbf{t}}-3s_{\mathbf{t}},\\[0,2mm]
0.5+\frac{g_{\mathbf{t}}^{z}(a)}{6}, & \mbox{if} & M_{\mathbf{t}}-3s_{\mathbf{t}}<g_{\mathbf{t}}(a)<M_{\mathbf{t}}+3s_{\mathbf{t}}, \\[0,2mm]
1, & \mbox{if} & g_{\mathbf{t}}(a)\geqslant M_{\mathbf{t}}+3s_{\mathbf{t}}
\end{array}
\right.
$$
\noindent if $g_{\mathbf{t}}$ has an increasing direction of preference (the more the evaluation $g_{\mathbf{t}}(a)$, the better $a$ is on $g_{\mathbf{t}}$), and as

$$
\overline{g}_{\mathbf{t}}(a)=
\left\{
\begin{array}{lll}
0 & \mbox{if} & g_{\mathbf{t}}(a)\geqslant M_{\mathbf{t}}+3s_{\mathbf{t}},\\[0,2mm]
0.5-\frac{g_{\mathbf{t}}^{z}(a)}{6}, & \mbox{if} & M_{\mathbf{t}}-3s_{\mathbf{t}}<g_{\mathbf{t}}(a)<M_{\mathbf{t}}+3s_{\mathbf{t}}, \\[0,2mm]
1, & \mbox{if} & g_{\mathbf{t}}(a)\leqslant M_{\mathbf{t}}-3s_{\mathbf{t}}
\end{array}
\right.
$$
\noindent if $g_{\mathbf{t}}$ has a decreasing direction of preference (the less the evaluation $g_{\mathbf{t}}(a)$, the better $a$ is on $g_{\mathbf{t}}$).

In consequence of the application of the considered normalization technique, the new data are therefore shown in Table \ref{NormalizationTable}. 

\begin{table}[!h]
  \centering
  \caption{Performances of the 28 Countries on the 11 elementary criteria taken into account}
	\resizebox{1\textwidth}{!}{
    \begin{tabular}{l|cccc|cccc|ccc}
		\hline
		\hline
          & \multicolumn{1}{c}{$g_{\mathbf{(1,1)}}$} & \multicolumn{1}{c}{$g_{\mathbf{(1,2)}}$} & \multicolumn{1}{c}{$g_{\mathbf{(1,3)}}$} & $g_{\mathbf{(1,4)}}$ & \multicolumn{1}{c}{$g_{\mathbf{(2,1)}}$} & \multicolumn{1}{c}{$g_{\mathbf{(2,2)}}$} & \multicolumn{1}{c}{$g_{\mathbf{(2,3)}}$} & $g_{\mathbf{(2,4)}}$ & $g_{\mathbf{(3,1)}}$ & $g_{\mathbf{(3,2)}}$ & $g_{\mathbf{(3,3)}}$ \\
		\hline			
		\hline
    Austria & 0.6176 & 0.6314 & 0.6144 & 0.4403 & 0.4970 & 0.6657 & 0.5290 & 0.4382 & 0.4487 & 0.5317 & 0.5224 \\
    Belgium & 0.5860 & 0.6229 & 0.5374 & 0.5708 & 0.5356 & 0.7126 & 0.4427 & 0.3362 & 0.5111 & 0.4293 & 0.5187 \\
    Bulgaria & 0.3062 & 0.4192 & 0.4900 & 0.4841 & 0.4954 & 0.2792 & 0.6404 & 0.6831 & 0.4644 & 0.5688 & 0.4671 \\    
    Cyprus & 0.4401 & 0.4192 & 0.1760 & 0.4845 & 0.8109 & 0.3551 & 0.3630 & 0.3624 & 0.4963 & 0.1157 & 0.1896 \\    
    Croatia & 0.3459 & 0.4574 & 0.5107 & 0.4287 & 0.6625 & 0.5548 & 0.3239 & 0.4324 & 0.4425 & 0.5611 & 0.2593 \\    
    Denmark & 0.6840 & 0.4828 & 0.6618 & 0.4407 & 0.3763 & 0.7298 & 0.6397 & 0.6477 & 0.4530 & 0.7489 & 0.6002 \\
    Estonia & 0.3955 & 0.6398 & 0.5789 & 0.5357 & 0.3872 & 0.4266 & 0.8115 & 0.7633 & 0.4863 & 0.4715 & 0.4494 \\
    Finland & 0.6024 & 0.5592 & 0.4545 & 0.3660 & 0.3330 & 0.7816 & 0.6699 & 0.6482 & 0.4169 & 0.3777 & 0.5216 \\
    France & 0.5456 & 0.5677 & 0.4633 & 0.3408 & 0.2000 & 0.8013 & 0.5652 & 0.3654 & 0.4088 & 0.3815 & 0.4966 \\
    Germany & 0.5971 & 0.4574 & 0.6440 & 0.4073 & 0.5913 & 0.5203 & 0.6124 & 0.4954 & 0.4338 & 0.7514 & 0.6370 \\    
    Greece & 0.3913 & 0.0499 & 0.1227 & 0.3462 & 0.6424 & 0.6707 & 0.4066 & 0.0378 & 0.3990 & 0.4101 & 0.3759 \\    
    Ireland & 0.7762 & 0.9921 & 0.8010 & 0.7483 & 0.6037 & 0.1061 & 0.2391 & 0.4878 & 0.6029 & 0.4365 & 1 \\    
    Italy & 0.4944 & 0.3300 & 0.4011 & 0.3458 & 0.6269 & 0.6337 & 0.2258 & 0.2284 & 0.4032 & 0.5304 & 0.5630 \\    
    Latvia & 0.3651 & 0.5168 & 0.4456 & 0.4626 & 0.4490 & 0.3477 & 0.5718 & 0.6168 & 0.4597 & 0.4033 & 0.3729 \\    
    Lithuania & 0.3700 & 0.3513 & 0.3271 & 0.5311 & 0.4444 & 0.2984 & 0.5637 & 0.6013 & 0.4807 & 0.3973 & 0.4331 \\    
    Luxembourg & 1     & 0.3810 & 0.5048 & 1 & 0.5016 & 0.4759 & 0.7717 & 0.8150 & 1     & 0.6234 & 0.4166 \\   
    Malta &  0.4569 & 0.5422 & 0.7625 & 0.7784 & 0.7831 & 0.3502 & 0.4560 & 0.5606 & 0.7434 & 0.8346 & 0.2617 \\    
    Netherlands & 0.6229 & 0.4743 & 0.6944 & 0.5683 & 0.6300 & 0.4759 & 0.6640 & 0.5223 & 0.5451 & 0.8052 & 0.7008 \\    
    Poland & 0.3469 & 0.4616 & 0.4011 & 0.4391 & 0.4057 & 0.4315 & 0.4523 & 0.5394 & 0.4403 & 0.3964 & 0.5112 \\    
    Portugal & 0.4085 & 0.2876 & 0.3212 & 0.3941 & 0.7367 & 0.5277 & 0.1764 & 0.2536 & 0.4293 & 0.4574 & 0.4346 \\    
    U.K.  & 0.5525 & 0.3470 & 0.1879 & 0.3371 & 0.3748 & 0.4167 & 0.3217 & 0.3952 & 0.4054 & 0.2031 & 0.4132 \\    
    Czech Rep. & 0.4010 & 0.7247 & 0.6618 & 0.5311 & 0.6161 & 0.3847 & 0.6677 & 0.6246 & 0.4798 & 0.5069 & 0.5984 \\    
    Romania & 0.3237 & 0.6398 & 0.4426 & 0.3866 & 0.1041 & 0.3329 & 0.5113 & 0.6181 & 0.4208 & 0.2787 & 0.4139 \\
    Slovakia & 0.3804 & 0.5804 & 0.4870 & 0.6075 & 0.3284 & 0.4512 & 0.5135 & 0.5528 & 0.5072 & 0.4016 & 0.5458 \\   
    Slovenia &  0.4272 & 0.4107 & 0.5433 & 0.5435 & 0.6021 & 0.5153 & 0.4080 & 0.4999 & 0.4855 & 0.6520 & 0.5675 \\   
    Spain &  0.4647 & 0.4871 & 0.4930 & 0.3577 & 0.2727 & 0.4315 & 0.3328 & 0.3719 & 0.4111 & 0.5023 & 0.4848 \\    
    Sweden & 0.6679 & 0.6823 & 0.7122 & 0.4011 & 0.5526 & 0.6337 & 0.7429 & 0.6429 & 0.4346 & 0.6204 & 0.5461 \\    
    Hungary & 0.3478 & 0.4828 & 0.5581 & 0.5827 & 0.4351 & 0.6879 & 0.3756 & 0.4581 & 0.5025 & 0.6012 & 0.5618 \\
		\hline
		\hline
    \end{tabular}%
		}
  \label{NormalizationTable}%
\end{table}%

For the sake of simplicity, we grouped the 25 rating classes to which Countries can be assigned by S\&P into four classes only, as shown in Table \ref{RatingCategories}. For example, a country that S\&P would assign to a class between AA- and AAA is assigned to $C_{4}$, while a country is assigned to category $C_{1}$ if its S\&P rating is lower than $BBB-$. In this way, Austria is assigned to $C_{4}$ being rated $AA+$, while Cyprus and Croatia are assigned to category $C_{1}$ since their ratings are $BB+$ and $BB$, respectively. The categories to which the considered Countries are assigned by using the described classification rule are shown in the last column of Table \ref{EvaluationsTable}. 

\begin{table}[!h]
  \centering
  \caption{Categories of rating defined for our case study}
	\resizebox{0.3\textwidth}{!}{
    \begin{tabular}{c|cc}
		\hline
		\hline
     Category Label & From... & To... \\
		\hline
		\hline
		$C_{4}$ & $AA-$ & $AAA$ \\
		$C_{3}$ & $A-$  & $A+$ \\
		$C_{2}$ & $BBB-$ & $BBB+$ \\
		$C_{1}$ & \multicolumn{2}{c}{otherwise} \\
 		\hline			
		\hline
    
		\hline
		\hline
    \end{tabular}%
		}
  \label{RatingCategories}%
\end{table}%

\subsection{Part I: Justification of the use of the Choquet integral for sorting}\label{FirstPart}
In financial applications, in general, a weighted sum is used as underlying preference model to aggregate the performances of the alternatives on the criteria at hand. Consequently, at first, we shall check if 
$$
\displaystyle WS(a,\mathbf{w})=\sum_{\mathbf{t}\in E_G}w_{\mathbf{t}}\overline{g}_{\mathbf{t}}(a)
$$
\noindent is able to represent the class assignments in Table \ref{EvaluationsTable}, where $\mathbf{w}=\left[w_{\mathbf{t}}\right]_{\mathbf{t}\in E_{G}}$ is such that $w_{\mathbf{t}}\geqslant 0$ for all $\mathbf{t}\in E_G$ and $\displaystyle\sum_{\mathbf{t}\in E_G}w_{\mathbf{t}}=1$. 

\noindent From a technical point of view, the existence of such a compatible model can be checked by solving the following LP problem

\begin{equation}\label{WSLPproblem}
\begin{array}{l}
\;\;\;\varepsilon^{*}_{WS}=\max\varepsilon, \;\;\mbox{subject to}, \\[2mm]
\left.
\begin{array}{l}
\left.
\begin{array}{l}
b_{h-1}\leqslant WS(a,\mathbf{w}),\\[2mm]
WS(a,\mathbf{w})+\varepsilon\leqslant b_{h}, \\[2mm]
\end{array}
\right\}\;\mbox{if $a\xrightarrow[{\mathbf{0}}]{} C_{h}$} \\[2mm]
\;w_{\mathbf{t}}\geqslant 0 \;\mbox{for all}\;\mathbf{t}\in E_G, \\[2mm]
\;\displaystyle\sum_{\mathbf{t}\in E_G}w_{\mathbf{t}}=1, \\[2mm]
\;b_{0}=0, \;b_{4}=1, \\[2mm]
\;b_{h-1}+\varepsilon\leqslant b_h, \;h=1,\ldots,4.
\end{array}
\right\}E^{DM}_{WS}
\end{array}
\end{equation}
\noindent We get that $E_{WS}^{DM}$ is feasible but $\varepsilon^{*}_{WS}=-0.0286<0$. This means that there not exists any vector of weights $\mathbf{w}$ and any threshold vector $\mathbf{b}=\left(b_0,b_1,\ldots,b_4\right)$, with $0=b_0<b_1<\ldots<b_4=1\footnote{Since all evaluations in Table \ref{NormalizationTable} are in the interval $[0,1]$ and, consequently, the maximum evaluation of a Country at comprehensive level is 1, we put $b_{4}=1$.}$, restoring the countries' classifications provided in Table \ref{EvaluationsTable}. Consequently, a weighted sum is not able to justify these assignments.

For this reason, we use a slightly more complex preference model, that is the 2-additive Choquet integral, described in Section \ref{MCHPChoquet}. Solving the LP problem 

\begin{equation}\label{ChLP}
\begin{array}{l}
\varepsilon^{*}_{Ch}=\max\varepsilon \;\mbox{subject to}, \\[3mm]
E_{Ch}^{DM}\cup E_{Ch}^{TC} 
\end{array}
\end{equation}

\noindent we get that the set $E_{Ch}^{DM}\cup E_{Ch}^{TC}$ is feasible and $\varepsilon_{Ch}^{*}=0.0024>0$. Consequently, there exists at least on set of M\"{o}bius parameters and a vector of thresholds $\mathbf{b}$ compatible with the considered Countries' assignments. 

The application of the 2-additive Choquet integral involves 66 parameters: 11 M\"{o}bius parameters for the elementary criteria and one M\"{o}bius parameter for each pair of elementary criteria. Therefore, as described in detail in Section \ref{MinimalSets}, we look for the minimal sets of interacting criteria necessary to restore the considered assignments. Solving the different MILP problems, we get the six different sets of pairs of interacting criteria shown in Table \ref{MinimalInteractingCriteria}. 

\begin{table}[htbp]
  \centering
  \caption{Minimal sets of pairs of interacting criteria. ``+" in the table represents a positive interaction, while ``-" represents a negative interactions.}
	\resizebox{0.6\textwidth}{!}{
    \begin{tabular}{lcccccc}
		\hline
		\hline
     & $\Gamma_{1}^{Min}$ & $\Gamma_{2}^{Min}$ & $\Gamma_{3}^{Min}$ & $\Gamma_{4}^{Min}$ & $\Gamma_{5}^{Min}$ & $\Gamma_{6}^{Min}$ \\
		\hline
		\hline
		$\{GDPc\;,\;Ep/GDP\}$  & + & + & + & + & + & + \\[1mm]
		$\{GDPc\;,\;I-Ex/R\}$  & - & - &   &   & - & - \\[1mm]
		$\{Ep/GDP\;,\;I-Ex/R\}$ & + & + & - & + &   &   \\[1mm]
		$\{Ep/GDP\;,\;TB/GDP\}$  & + &   &   &   & + & + \\[1mm]
		$\{CAR/GDP\;,\;TB/GDP\}$    &   & + & + & + &   &   \\[1mm]
		$\{I-Ex/R\;,\;TB/GDP\}$  &   &   & + &   &   &   \\[1mm]
		$\{GDPc\;,\;D/GDP\}$ &   &   &   & - & + &   \\[1mm]
		$\{D/GDP\;,\;CAR/GDP\}$   &   &   &   &   &   & + \\ 
		\hline			
		\hline
    \end{tabular}%
		}
  \label{MinimalInteractingCriteria}%
\end{table}%

Looking at the table, one can observe that 4 pairs of interacting criteria are necessary to justify the considered assignments. This is a relevant reduction since, in this way, only 15 M\"{o}bius parameters are enough in the considered preference model instead of the 66 mentioned above. Moreover, many times the same pair of interacting criteria is present in all discovered minimal sets and, above all, the \textit{core} of these minimal sets of interacting criteria \citep{MUSA_int,GrecoMatarazzoSlowinski2001} is represented by the pair $\{GDPpc,ExGDP\}$ that is composed of positively interacting criteria in all minimal sets. This means that the positive interaction between GDP per capita and Exports/GDP is fundamental to explain the provided Countries' assignments. \\
At the same time one can underline that, as observed in Section \ref{MinimalSets}, some pairs of criteria are positively interacting in one minimal set and negatively interacting in another. For example, $Ep/GDP$ and $I-Ex/R$ are negatively interacting in $\Gamma_{3}^{Min}$ and positively interacting in the other three minimal sets in which they interact, while, $GDPc$ and $D/GDP$ are negatively interacting in $\Gamma_{4}^{Min}$ and positively interacting in $\Gamma_{5}^{Min}$.

\subsection{Part II: Exploring the results of NAROR and SMAA}\label{SecondPart}%
While, in the previous section, we have shown the potentialities of the 2-additive Choquet integral to represent the assignments provided by the DM, in this section we shall show the application of the NAROR and SMAA methodologies to sorting problems when the 2-additive Choquet integral is considered as preference model. 

Differently from the case above in which the preference assignments regarded all the 28 EU countries, in this section we shall assume that the DM provided the following preference information in terms of:
\begin{itemize}
\item comparison between criteria:
\begin{itemize}
\item Net debt/GDP is more important than Primary balance/GDP,
\item Savings/GDP is more important than Investment/GDP,
\item Current account balance/GDP is more important than Current account receipts/GDP,
\end{itemize}

\item interaction between criteria: 
\begin{itemize}
\item There is a negative interaction between GDP per capita and Exports/GDP,
\item There is a positive interaction between Interest/revenues and Net Debt/GDP,
\item There is a positive interaction between Current account receipts/GDP and Current account balance/GDP,
\end{itemize}

\item assignments examples:
\begin{table}[!htb]
\smallskip
\centering
\begin{small}
\resizebox{0.6\textwidth}{!}{
\begin{tabular}{ccrrrrcc} 
		\hline
		\hline   
    \multicolumn{2}{c}{\textit{Global}} & \multicolumn{2}{c}{\textit{Ec}} & \multicolumn{2}{c}{\textit{Gov}} & \multicolumn{2}{c}{\textit{Fin}} \\
		\hline
		\hline
    Germany & $C_4$    & \multicolumn{1}{c}{Sweden} & \multicolumn{1}{c}{$C_4$} & \multicolumn{1}{c}{Netherlands} & \multicolumn{1}{c}{$C_4$} & Denmark & $C_4$ \\
    Slovenia & $C_3$    & \multicolumn{1}{c}{Spain} & \multicolumn{1}{c}{$C_2$} & \multicolumn{1}{c}{Slovakia} & \multicolumn{1}{c}{$C_3$} & Lithuania & $C_3$ \\
    Italy & $C_2$    & \multicolumn{1}{c}{Cyprus} & \multicolumn{1}{c}{$C_1$} & \multicolumn{1}{c}{Portugal} & \multicolumn{1}{c}{$C_2$} & Poland & $C_2$ \\
    Greece  & $C_1$    &       &       &       &       & Greece & $C_1$ \\
		\hline
		\hline
    \end{tabular}%
}
\end{small}
\end{table}
\end{itemize}

\subsubsection{NAROR results}\label{narorresult}
Solving the LP problem (\ref{ChLP}), we get $\varepsilon^{*}_{Ch}>0$. This means that there exists at least one set of M\"{o}bius parameters and at least one thresholds vector compatible with the preference information provided in the previous section. To get more stable assignments, in this section, we shall apply the NAROR. In Tables \ref{necessglob} and \ref{possglob} we presented the \textit{new} necessary and possible assignments of the countries on which the DM did not provide any assignment. 

\begin{table}[!htb]
\smallskip
\centering
\caption{New necessary assignments at global level as well as with respect to Economic, Governmental and Financial points of view\label{necessglob}}
\begin{small}
\resizebox{0.4\textwidth}{!}{
\begin{tabular}{ccccrrcc}
		\hline
		\hline
    \multicolumn{2}{c}{\textit{Global}} & \multicolumn{2}{c}{\textit{Ec}} & \multicolumn{2}{c}{\textit{Gov}} & \multicolumn{2}{c}{\textit{Fin}} \\
		\hline
		\hline
     &    & {Ireland} & {$C_4$}  &  &  & {Netherlands} & {$C_4$} \\
     &    & Greece & $C_1$ &  &  &  &  \\
		\hline
		\hline
    \end{tabular}%
}
\end{small}
\end{table}

\begin{table}[!htb]
\smallskip
\centering
\caption{New possible assignments at global level as well as with respect to Economic, Governmental and Financial points of view \label{possglob}}
\begin{small}
\resizebox{0.7\textwidth}{!}{
\begin{tabular}{cccccccc}
		\hline
		\hline
    \multicolumn{2}{c}{\textit{Global}} & \multicolumn{2}{c}{\textit{Ec}} & \multicolumn{2}{c}{\textit{Gov}} & \multicolumn{2}{c}{\textit{Fin}} \\
		\hline
		\hline
     & & Austria & $\left[C_{2},C_{4}\right]$ & Croatia    & $\left[C_{2},C_{4}\right]$ & Austria & $\left[C_{2},C_{4}\right]$ \\
     & & Belgium & $\left[C_{2},C_{4}\right]$ & Denmark    & $\left[C_{3},C_{4}\right]$ & Belgium & $\left[C_{3},C_{4}\right]$ \\
     & & Greece  & $C_{1}$                    & Finland    & $\left[C_{3},C_{4}\right]$ & Estonia & $\left[C_{3},C_{4}\right]$ \\
     & & Ireland & $C_{4}$                    & Luxembourg & $\left[C_{3},C_{4}\right]$ & Finland & $\left[C_{1},C_{2}\right]$ \\
     & &         &                            & Spain      & $\left[C_{1},C_{3}\right]$ & France  & $\left[C_{1},C_{2}\right]$ \\
     & &         &                            & Sweden     & $\left[C_{3},C_{4}\right]$ & Ireland & $\left[C_{3},C_{4}\right]$ \\
     & &         &                            &            &                            & Luxembourg & $\left[C_{3},C_{4}\right]$ \\
     & &         &                            &            &                            & Netherlands & $C_{4}$ \\
     & &         &                            &            &                            & U.K.        & $\left[C_{1},C_{2}\right]$ \\
     & &         &                            &            &                            & Czech Rep.  & $\left[C_{3},C_{4}\right]$ \\
     & &         &                            &            &                            & Romania & $\left[C_{1},C_{2}\right]$ \\
     & &         &                            &            &                            & Slovakia & $\left[C_{3},C_{4}\right]$ \\
     & &         &                            &            &                            & Slovenia & $\left[C_{3},C_{4}\right]$ \\
     & &         &                            &            &                            & Hungary & $\left[C_{3},C_{4}\right]$ \\
		\hline
		\hline
    \end{tabular}%
}
\end{small}
\end{table}

Looking at table \ref{necessglob} one can observe that the new necessary assignments obtained by the application of the Choquet integral are very few. Indeed, at global level and with respect to governmental aspects, there is not any new necessary assignments; with respect to economic aspects, Ireland is necessarily assigned to $C_{4}$, while Greece is assigned to $C_{1}$. Finally, on financial aspects, Netherlands has to be necessarily assigned to $C_{4}$. \\
Analogously, in Table \ref{possglob} we reported the new possible assignments regarding the countries which have not been provided as reference assignments from the DM. Again, the column corresponding to the global level is empty meaning that all countries on which the DM did not provide any preference at global level, can be possibly assigned to the whole range of classes, that is $\left[C_{1},C_{4}\right]$. Looking at the particular aspects, some more in deep information can be obtained. For example, on economic aspects, Austria can be possibly assigned to the interval of classes $\left[C_2,C_4\right]$. This means that there is not any instance of the preference model compatible with the preferences provided by the DM assigning Austria to the class $C_{1}$. On governmental aspects, 6 countries can be assigned to less than the four considered classes and, in particular, Croatia and Spain to three classes, while Denmark, Finland, Luxembourg and Sweden to the interval of classes $\left[C_{3},C_{4}\right]$. Finally, the financial aspect is the one on which more new information is obtained. Indeed, of the 24 countries on which the DM did not provide any preference information, 14 can be possibly assigned to an interval of classes narrower than $\left[C_{1},C_{4}\right]$. In particular, Austria can be assigned to three classes in the interval $\left[C_{2},C_{4}\right]$, while all the other countries, apart from the Netherlands that could be necessarily assigned (and consequently possibly assigned) to $C_{4}$, can be placed in two contiguous classes only. 

\subsubsection{SMAA results}\label{smaaresult}
Since most of the countries can be possibly assigned to more than one class at global level as well at partial level, in this section we shall describe the application of the SMAA methodology, as described in section \ref{SMAAsubsec}, which permitted to get more robust results in terms of frequency of assignments of any Country at global level as well as with respect to the 3 considered aspects in the considered hierarchy of criteria. In Table \ref{CAI}, we report the class acceptability indices of all the considered Countries. 

\begin{table}[!htb]
\smallskip
\centering
\caption{Class acceptability indices at global level as well as with respect to the three macro-criteria at hand \label{CAI}}
\resizebox{0.9\textwidth}{!}{
\begin{tabular}{l|cccc|cccc|cccc|cccc}
          			\hline
					\hline
          & \multicolumn{4}{c|}{\textbf{Global level}} & \multicolumn{4}{c|}{\textbf{Ec}} & \multicolumn{4}{c|}{\textbf{Gov}} & \multicolumn{4}{c}{\textbf{Fin}} \\ 
          & $C_{1}$ & $C_{2}$ & $C_{3}$ & $C_{4}$ & $C_{1}$ & $C_{2}$ & $C_{3}$ & $C_{4}$ & $C_{1}$ & $C_{2}$ & $C_{3}$ & $C_{4}$ & $C_{1}$ & $C_{2}$ & $C_{3}$ & $C_{4}$ \\
          \hline
          \hline
          Austria & 0     & 4.404 & 95.596 & 0     & 0     & 0     & 4.461 & 95.539 & 0     & 4.902 & 80.623 & 14.475 & 0     & 0     & 77.839 & 22.161 \\
    Belgium & 0     & 70.220 & 29.780 & 0     & 0     & 29.949 & 33.432 & 36.619 & 0     & 40.468 & 59.450 & 0.082 & 0     & 0     & 78.712 & 21.288 \\
    Bulgaria & 0     & 34.678 & 65.322 & 0     & 18.495 & 81.505 & 0     & 0     & 0     & 0     & 0.506 & 99.494 & 0     & 0     & 63.929 & 36.071 \\
    Cyprus & 100   & 0     & 0     & 0     & 100   & 0     & 0     & 0     & 0     & 0.207 & 75.155 & 24.638 & 100   & 0     & 0     & 0 \\
    Croatia & 0     & 87.318 & 12.682 & 0     & 3.348 & 93.771 & 2.881 & 0     & 0     & 0     & 75.305 & 24.695 & 0     & 0     & 100   & 0 \\
    Denmark & 0     & 0     & 0     & 100   & 0     & 0     & 0.562 & 99.438 & 0     & 0     & 0.167 & 99.833 & 0     & 0     & 0     & 100 \\
    Estonia & 0     & 1.625 & 98.375 & 0     & 0     & 15.083 & 27.477 & 57.440 & 0     & 0     & 0     & 100   & 0     & 0     & 82.225 & 17.775 \\
    Finland & 0     & 100   & 0     & 0     & 2.082 & 96.484 & 1.434 & 0     & 0     & 0     & 0.244 & 99.756 & 85.270 & 14.730 & 0     & 0 \\
    France & 76.851 & 23.149 & 0     & 0     & 1.490 & 97.241 & 1.269 & 0     & 0     & 99.561 & 0.439 & 0     & 93.813 & 6.187 & 0     & 0 \\
    Germany & 0     & 0     & 0     & 100   & 0     & 0.340 & 11.084 & 88.575 & 0     & 0     & 27.471 & 72.529 & 0     & 0     & 5.396 & 94.604 \\
    Greece & 100   & 0     & 0     & 0     & 100   & 0     & 0     & 0     & 0     & 100   & 0     & 0     & 100   & 0     & 0     & 0 \\
    Ireland & 0     & 0     & 4.750 & 95.250 & 0     & 0     & 0     & 100   & 0     & 23.676 & 75.045 & 1.279 & 0     & 0     & 7.050 & 92.950 \\
    Italy & 0     & 100   & 0     & 0     & 64.928 & 35.072 & 0     & 0     & 0     & 100   & 0     & 0     & 0     & 0     & 100   & 0 \\
    Latvia & 7.259 & 92.741 & 0     & 0     & 30.704 & 69.296 & 0     & 0     & 0     & 0     & 22.890 & 77.110 & 0.035 & 90.684 & 9.281 & 0 \\
    Lithuania & 15.579 & 84.421 & 0     & 0     & 100   & 0     & 0     & 0     & 0     & 0     & 45.113 & 54.887 & 0     & 0     & 100   & 0 \\
    Luxembourg & 0     & 0     & 0     & 100   & 0     & 32.017 & 27.526 & 40.457 & 0     & 0     & 0     & 100   & 0     & 0     & 0     & 100 \\
    Malta & 0     & 0     & 0     & 100   & 0     & 0     & 0     & 100   & 0     & 0     & 0     & 100   & 0     & 0     & 0     & 100 \\
    Netherlands & 0     & 0     & 0     & 100   & 0     & 0     & 0     & 100   & 0     & 0     & 0     & 100   & 0     & 0     & 0     & 100 \\
    Poland & 66.584 & 33.416 & 0     & 0     & 59.212 & 40.788 & 0     & 0     & 0     & 0     & 94.320 & 5.680 & 0     & 100   & 0     & 0 \\
    Portugal & 49.364 & 50.636 & 0     & 0     & 100   & 0     & 0     & 0     & 0     & 100   & 0     & 0     & 0     & 47.882 & 52.118 & 0 \\
    U.K.  & 100   & 0     & 0     & 0     & 100   & 0     & 0     & 0     & 0     & 100   & 0     & 0     & 100   & 0     & 0     & 0 \\
    Czech Rep. & 0     & 0     & 62.111 & 37.889 & 0     & 0     & 0     & 100   & 0     & 0     & 0     & 100   & 0     & 0     & 61.996 & 38.004 \\
    Romania & 99.478 & 0.522 & 0     & 0     & 18.087 & 81.895 & 0.018 & 0     & 0     & 62.665 & 37.335 & 0     & 100   & 0     & 0     & 0 \\
    Slovakia & 0     & 96.311 & 3.689 & 0     & 0     & 78.798 & 20.208 & 0.994 & 0     & 0.000 & 100.000 & 0     & 0     & 0     & 91.353 & 8.647 \\
    Slovenia & 0     & 0     & 100   & 0     & 0     & 63.647 & 27.815 & 8.537 & 0     & 0     & 50.404 & 49.596 & 0     & 0     & 20.577 & 79.423 \\
    Spain & 18.253 & 81.747 & 0     & 0     & 0     & 100   & 0     & 0     & 0     & 100   & 0     & 0     & 0     & 5.885 & 94.115 & 0 \\
    Sweden & 0     & 0     & 0     & 100   & 0     & 0     & 0     & 100   & 0     & 0     & 0     & 100   & 0     & 0     & 54.529 & 45.471 \\
    Hungary & 0     & 10.564 & 89.436 & 0     & 0     & 44.028 & 33.960 & 22.012 & 0     & 19.906 & 80.094 & 0     & 0     & 0     & 31.537 & 68.463 \\
    \hline
    \hline
    \end{tabular}
}
\end{table}

Looking at the results in the table, one can observe that the number of classes to which each Country can be assigned is narrower than the interval of classes to which the same Countries can be possibly assigned by the application of the NAROR. For example, we observed that Austria can be possibly assigned to the interval $\left[C_{2},C_{4}\right]$ at global level. Anyway, looking at the class acceptability indices in Table \ref{CAI} Austria is assigned to two classes only and, in particular, to $C_{2}$ and $C_{3}$ with frequencies $4.404\%$ and $95.596\%$, respectively. This is not surprising since, as explained in the section devoted to the SMAA methodology, the CAI are obtained performing a sampling of the instances of the model compatible with the preferences provided by the DM. Consequently, the fact that Austria is assigned to $C_{4}$ with a frequency equal to $0\%$ means that among the instances sampled by SMAA there is not anyone assigning Austria to this class. Analogous observations can be done at macro-criteria level observing that the number of classes to which each country can be assigned is never 4 and most of the countries have a class acceptability index greater than 0\% for two classes only. 

Since, as previously explained, for each instance of the considered preference model compatible with the information provided by the DM we can get an assignment of the countries at global level as well as at macro-criteria level, we have computed the most frequent classification of the countries, that is the classification that appeared more frequently applying all the compatible instances. We got 124 different countries' assignments at comprehensive level, 437 countries' assignments on economic aspects, 205 countries' assignments on economic aspects and, finally, 91 countries' assignments on financial aspects. The most frequent ones are reported in Table \ref{MostFrequentClassifications}.

\begin{table}[!htb]
\smallskip
\centering
\caption{The most frequent assignments at global level as well as with respect to the three macro-criteria taken into account. The frequencies are referring to the number of times the considered assignment was obtained considering the whole set of compatible models.\label{MostFrequentClassifications}} 
\begin{small}
\resizebox{0.5\textwidth}{!}{
  \begin{tabular}{lcccc}
	\hline
	\hline
    & $Global$ & ${Ec}$ & ${Gov}$ & ${Fin}$ \\
	\hline
	\hline
	\multicolumn{1}{l}{Austria} & $C_3$     & $C_4$     & $C_3$     & $C_3$ \\
    \multicolumn{1}{l}{Belgium} & $C_2$     & $C_2$     & $C_3$     & $C_3$ \\
    \multicolumn{1}{l}{Bulgaria} & $C_3$     & $C_2$     & $C_4$     & $C_3$ \\
    \multicolumn{1}{l}{Cyprus} & $C_{1}$     & $C_{1}$     & $C_3$     & $C_{1}$ \\
    \multicolumn{1}{l}{Croatia} & $C_2$     & $C_2$     & $C_3$     & $C_3$ \\
    \multicolumn{1}{l}{Denmark} & $C_4$     & $C_4$     & $C_4$     & $C_4$ \\
    \multicolumn{1}{l}{Estonia} & $C_3$     & $C_2$     & $C_4$     & $C_3$ \\
    \multicolumn{1}{l}{Finland} & $C_2$     & $C_2$     & $C_4$     & $C_1$ \\
    \multicolumn{1}{l}{France} & $C_{1}$     & $C_2$     & $C_2$     & $C_{1}$ \\
    \multicolumn{1}{l}{Germany} & $C_4$     & $C_4$     & $C_3$     & $C_4$ \\
    \multicolumn{1}{l}{Greece} & $C_{1}$     & $C_{1}$     & $C_2$     & $C_{1}$ \\
    \multicolumn{1}{l}{Ireland} & $C_4$     & $C_4$     & $C_3$     & $C_4$ \\
    \multicolumn{1}{l}{Italy} & $C_2$     & $C_2$     & $C_2$     & $C_3$ \\
    \multicolumn{1}{l}{Latvia} & $C_2$     & $C_2$     & $C_3$     & $C_2$ \\
    \multicolumn{1}{l}{Lithuania} & $C_2$     & $C_{1}$     & $C_3$     & $C_3$ \\
    \multicolumn{1}{l}{Luxembourg} & $C_4$     & $C_2$     & $C_4$     & $C_4$ \\
    \multicolumn{1}{l}{Malta} & $C_4$     & $C_4$     & $C_4$     & $C_4$ \\
    \multicolumn{1}{l}{Netherlands} & $C_4$     & $C_4$     & $C_4$     & $C_4$ \\
    \multicolumn{1}{l}{Poland} & $C_{1}$     & $C_2$     & $C_3$     & $C_2$ \\
    \multicolumn{1}{l}{Portugal} & $C_{1}$     & $C_{1}$     & $C_2$     & $C_2$ \\
    \multicolumn{1}{l}{U.K.} & $C_{1}$     & $C_{1}$     & $C_2$     & $C_{1}$ \\
    \multicolumn{1}{l}{Czech Rep.} & $C_4$     & $C_4$     & $C_4$     & $C_3$ \\
    \multicolumn{1}{l}{Romania} & $C_{1}$     & $C_2$     & $C_3$     & $C_1$ \\
    \multicolumn{1}{l}{Slovakia} & $C_2$     & $C_2$     & $C_3$     & $C_3$ \\
    \multicolumn{1}{l}{Slovenia} & $C_3$     & $C_2$     & $C_3$     & $C_4$ \\
    \multicolumn{1}{l}{Spain} & $C_2$     & $C_2$     & $C_2$     & $C_3$ \\
    \multicolumn{1}{l}{Sweden} & $C_4$     & $C_4$     & $C_4$     & $C_3$ \\
    \multicolumn{1}{l}{Hungary} & $C_3$     & $C_2$     & $C_3$     & $C_4$ \\
    \hline
    \hline
    Frequency & 9.904\% & 10.416\% & 11.218\%  & 10.607 \\
    \hline
    \hline
    \end{tabular}%
}
\end{small}
\end{table}

As already described in Section \ref{SMAAsubsec}, the class acceptability indices shown in Table \ref{CAI} can be aggregate to produce a single precise assignment by minimizing the function in eq. (\ref{LossFunc}). In particular, let us show that the final assignment is dependent on the choice of the distance function $d$ considering the economic macro-criterion. For such a reason, let us take into account three different distances being $d\left(C_h,C_k\right)=1$, $d\left(C_h,C_k\right)=|h-k|$ and $d\left(C_h,C_k\right)=\sqrt{|h-k|}$. On one hand, the choice $d\left(C_h,C_k\right)=1$ implies that the error in assigning an alternative to a wrong class is not dependent on how far is the class to which the alternative is wrongly assigned; on the other hand, $d\left(C_h,C_k\right)=|h-k|$ and $d\left(C_h,C_k\right)=\sqrt{|h-k|}$ imply that the error increases more in case $d\left(C_h,C_k\right)=|h-k|$ is considered and less in considering $d\left(C_h,C_k\right)=\sqrt{|h-k|}$. Minimizing the function in eq. (\ref{LossFunc}) we get that, independently on the choice of the three distances $d$, all Countries are assigned to the class $h$ for which the corresponding class acceptability index is the highest, apart from the three Countries shown in Table \ref{texlossddd}. Indeed, in the case $d\left(C_h,C_k\right)=1$, the three countries are assigned to the class in correspondence of which the Country presents the highest CAI, so Belgium and Luxembourg are assigned to $C_4$, while Hungary is assigned to $C_2$. In the other two cases, the Countries are assigned to different classes. In particular, Belgium and Luxembourg are assigned in both cases to $C_{3}$, while Hungary is assigned to $C_3$, when $d\left(C_h,C_k\right)=|h-k|$ and again to $C_{2}$, when $d\left(C_h,C_k\right)=\sqrt{|h-k|}.$

\begin{small}
\begin{table}[h]
\centering
  \caption{Final assignments for three different distances $d$\label{texlossddd}}
	  \begin{footnotesize}
		\begin{tabular}{lccc}
		\hline
		\hline
		 & $d=1$ & $d=|h-k|$ & $\sqrt{|h-k|}$ \\
		\hline
		\hline
		Belgium     &  ${C}_{4}$   &  ${C}_{3}$  &  ${C}_{3}$ \\
		Luxembourg  &  ${C}_{4}$   &  ${C}_{3}$  &  ${C}_{3}$ \\
		Hungary     &  ${C}_{2}$   &  ${C}_{3}$  &  ${C}_{2}$ \\
    \hline
		\hline
	  \end{tabular}
		\end{footnotesize}
\end{table}
\end{small}

Let us conclude this section by observing one thing that can be also the introduction for the next section. Of course, the results obtained by the NAROR and the SMAA methodologies are strictly dependent on the exemplary assignments that have been provided by the DM. For example, looking at the most frequent country assignments shown in Table \ref{MostFrequentClassifications}, no country has been assigned to $C_{3}$ on economic aspects as well as no country has been assigned to $C_{1}$ on governmental aspects. In particular, this result is also enforced by the fact that the class assignment index of all countries on $C_{1}$ at governmental level is always equal to $0\%$. Indeed, looking at the preferences provided by the DM, one can observe that the DM was not able to assign any country to the class $C_{3}$ on economic aspects as well as any country to $C_{1}$ on governmental aspects. So the question is now ``which is the predictive capability of the method"? In other words, we are asking which is the capability of the method to get a correct countries' assignment starting from the example assignments provided by the DM. This aspect will be investigated in the next section.

\subsection{Part IV: Predictive capability of the proposed method}\label{ThirdPart}
As already introduced in the section above, in this section we shall study the predictive capability of the proposed sorting model as a function of the number of assignments provided by the DM at comprehensive level. To this aim we propose a comparison of four methods by means of an \textit{in-sample forecast} procedure \citep{brooks2014introductory} in which the parameters of the model are estimated by some examples and then applied to the whole set. The reasons to adopt this procedure are mainly two: at first, the results of the procedure are easy to be read and they permit a fair comparison between the considered methods. At second, the in-sample forecast procedures have shown a significant predictive ability in the financial framework \citep{inoue2005sample,rapach2006sample}.

We have compared the following four methods: 
\begin{itemize}
\item weighted sum ($WS$),
\item 2-additive Choquet integral with all possible pairs of interacting criteria ($CH$),
\item 2-additive Choquet integral with minimal sets of pairs of interacting criteria ($MSCH$),
\item Multiple linear regression model (MLR) \citep{montgomery2012introduction}.
\end{itemize}

With respect to MLR we consider a multiple linear regression model similar to the one presented in \cite{ferri1999procyclical} in which we use as explanatory variables the eleven elementary criteria introduced in Section \ref{CaseStudy}. So, we suppose that the assignments have to be explained by the following model:

\begin{footnotesize}
\begin{equation}\label{MLRmodel}
{S}_{j} = \alpha + {\beta}_{1}{GDPc_{i}} + {\beta}_{2}{\frac{(I/GDP)_{i}}{GDPc_{i}}} +{\beta}_{3} {\frac{(S/GDP)_{i}}{GDPc_{i}}} +{\beta}_{4}{\frac{(Ep/GDP)_{i}}{GDPc_{i}}} +{\beta}_{5}{\frac{(PB/GDP)_{i}}{GDPc_{i}}} +
\end{equation}
$$
+{\beta}_{6}{\frac{(Ex/GDP)_{i}}{GDPc_{i}}}+{\beta}_{7}{\frac{(I-Ex/R)_{i}}{GDPc_{i}}} +{\beta}_{8}{\frac{(D/GDP)_{i}}{GDPc_{i}}}+{\beta}_{9}{\frac{(CAR/GDP)_{i}}{GDPc_{i}}} +{\beta}_{10}{\frac{(CAB/GDP)_{i}}{GDPc_{i}}} +{\beta}_{11}{\frac{(TB/GDP)_{i}}{GDPc_{i}}} +{\varepsilon}
$$
\end{footnotesize}

\noindent where $j=1,2,3,4$ denotes the class to which the country $i=1,\ldots,28$, can be assigned and index $i$ is used to denote the Country at hand. For example, $GDPc_{i}$ is the evaluation of country $i$ on elementary criterion $GDPc$. 

Regarding the MLR model, the following procedure is performed: 
\begin{description}
\item[Step 1] we sample $k$ reference assignments at comprehensive level of the 28 countries at hand shown in Table \ref{EvaluationsTable}, with $k=1,\ldots,27$,
\item[Step 2] using the $k$ reference assignments, we estimate all the parameters of the model (\ref{MLRmodel}), that is $\alpha,$ $\beta_{t}$, $t=1,\ldots,11$, and $\varepsilon$, 
\item[Step 3] we apply the considered model using the inferred parameters to assign the remaining $28-k$ countries. The assignment is performed as follows: The assignment of the predicted value is done in terms of minimum distance (absolute value) from one of the four classes.
\end{description}
The procedure described above is therefore performed 50 times for each $k=1,\ldots,27$ and 28 times for $k=27$. This means that 50 times we sampled $k$ assignments, different in each sample so that a single sample of $k$ assignments is not considered more than once. In the case $k=27,$ the procedure is not a sampling since only 1 country has not been assigned and, therefore, it can be performed 28 times (once for each country).

 To evaluate the goodness of the fit, we have considered the coefficients of determination \citep{kvaalseth1985cautionary,seber2012linear} (adjusted ${R}^{2}$ or ${\overline{R}}^{2}$) and the \textit{p}-value of the \textit{F}-statistic \citep{seber2012linear} (at the 5\% significance level). The results of the regressions validation are shown in Table \ref{validation}, where for each number of reference assignments, we report the minimum, medium and maximum value of ${\overline{R}}^{2}$ and the minimum, medium and maximum \textit{p}-value of the \textit{F}-statistic.

\begin{table}[!htb]
\centering
\caption{Values of the regressions validation\label{validation}}
\resizebox{0.7\textwidth}{!}{
\begin{tabular}{ccc}
		\hline
		\hline
    $\#$ reference assignments & ${\overline{R}}^{2}$ (min/med/max) & \textit{p}-value of F-stat. (min/med/max) \\
    \hline
		\hline
    14    & -0.1954 / 0.6865 / 0.9894  &     0.0089 / 0.2179 / 0.6730 \\
    16    & -0.0036 / 0.6247 / 0.9065  &     0.0103 / 0.1476 / 0.5543 \\
    18    &  0.2405 / 0.6115 / 0.9220  & 8.4956e-04 / 0.0915 / 0.3245 \\
    20    &  0.2719 / 0.5766 / 0.8722  & 6.5330e-04 / 0.0731 / 0.2451 \\
    22    &  0.3523 / 0.5716 / 0.8477  & 2.7331e-04 / 0.0429 / 0.1360 \\
    24    &  0.4107 / 0.5440 / 0.7643  & 6.6228e-04 / 0.0276 / 0.0688 \\
    26    &  0.4688 / 0.5573 / 0.7609  & 2.2877e-04 / 0.0128 / 0.0280 \\
    27    &  0.4914 / 0.5408 / 0.7069  & 5.1416e-04 / 0.0104 / 0.0174 \\
		\hline
		\hline
    \end{tabular}%
}
\end{table}

With respect to the other three models, the procedure is a bit different. Indeed, while for the MLR we are estimating a single vector of values, that is a single value for each of the considered parameters in (\ref{MLRmodel}), for the other three models we are implementing a procedure aiming to take into account robustness concerns and composed, therefore, of the following steps:

\begin{description}
\item[Step 1'] we sample $k$ reference assignments at comprehensive level of the 28 countries at hand shown in Table \ref{EvaluationsTable}, with $k=1,\ldots,27$,
\item[Step 2'] we check, for each of the three methods, if there exists at least one instance of the preference model compatible with these reference assignments. Let us observe that if $WS$ is able to restore the reference assignments then, of course, $CH$ is able to restore the same assignments too since the $WS$ is a particular case of $CH$ when there is not any interaction between the considered criteria. Moreover, in this case we will not take into account $MSCH$ since the minimal number of pairs of interacting criteria is equal to zero. For this reason, if $WS$ is able to restore the preferences provided by the DM, then we shall consider models $WS$ and $CH$ in the following steps, while, in the opposite case ($WS$ is not able to restore the preferences provided by the DM), in the following steps we shall consider models $CH$ and $MSCH$. In particular, regarding $MSCH$, we shall check all minimal sets of pairs of interacting criteria compatible with the provided reference assignments. \\
If neither $WS$ nor $CH$ are able to restore the reference assignments and we have not performed yet the 50 runs\footnote{28 in the case $k=27$}, we go again to Step 1' sampling another set of $k$ reference assignments, otherwise we go to Step 3',
\item[Step 3'] for each of the two considered models compatible with the preferences provided by the DM ($WS$ and $CH$ or $CH$ and $MSCH$) we sample 10,000 compatible instances and we perform the assignment of all countries with respect to all sampled compatible models. Regarding $MSCH$, we sample 10,000 compatible models for each considered minimal set of pairs of interacting criteria. This means that if 3 minimal sets are able to restore the reference assignments, then we shall sample 30,000 compatible instances getting, therefore, 30,000 different countries' assignments, 
\item[Step 4'] for each of the two considered models compatible with the assignments provided by the DM ($WS$ and $CH$ or $CH$ and $MSCH$) we compute the CAI of all countries at hand considering the assignments performed at the previous step,
\item[Step 5'] each country $a$ is considered correctly assigned if $CAI(a,h_{correct})\geqslant p\%$, where $p\in\{50,75,100\}$ and $C_{h_{correct}}$ is the class to which country $a$ had to be assigned in Table \ref{EvaluationsTable}. We are therefore assuming that the model is assigning correctly Country $a$ if the CAI corresponding to the class to which the DM had liked to assign it is no lower than $p\%$. 
\end{description}

\begin{figure}[htb!]
\centering
\caption{Percentage of correct assignments for the non-reference countries considering each of the four methods}
\label{ComparisonFigure}
\includegraphics[scale=0.35]{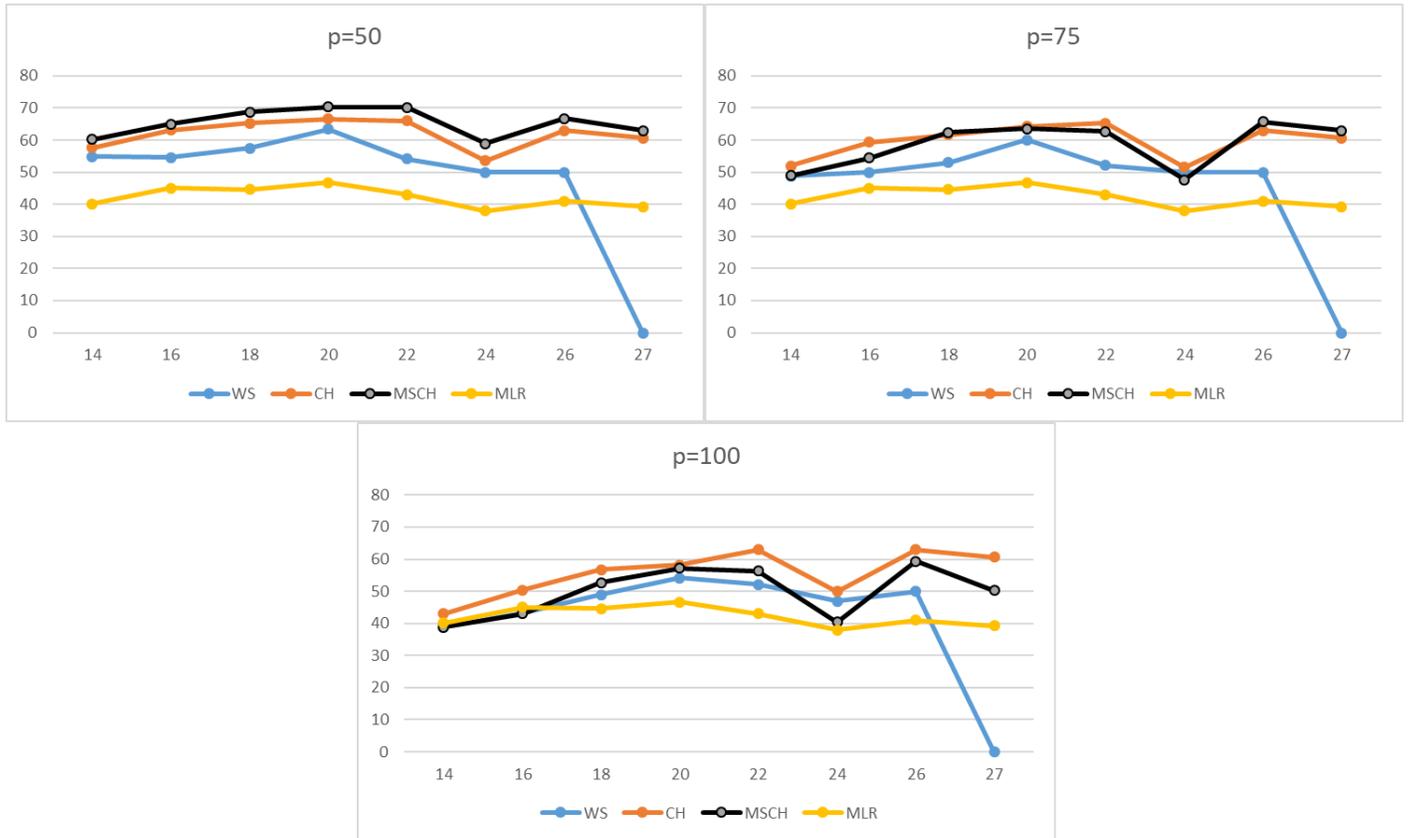}
\end{figure}

In Table \ref{ComparisonTable} and in Fig. \ref{ComparisonFigure}, we report and show the percentage of non-reference countries that have been assigned to the correct class by means of each of the four methods at hand and considering the three values of $p$, that is $p=50,$ $p=75$ and $p=100$. 

Before commenting more in depth the results, let us underline from the data in Table \ref{CompatibleModels} that $CH$ was able in all sampling to restore the considered reference assignments, while it was not the case for the $WS$. In particular, the number of runs out of the 50 considered in which $WS$ was able to restore the reference assignments is decreasing with respect to $k$ passing from its maximum $(43)$ when $k=14$ to its minimum $(1)$ when $k=27$. Another consideration is that while all the percentages in Table \ref{ComparisonTable} (and, consequently, in Fig. \ref{ComparisonFigure}) are obtained by the average over the 50 different samplings, the data for $WS$ and $MSCH$ are obtained by the average of the number or times the $WS$ was able or not to restore the reference assignments. For example, while the percentages corresponding to the case $k=14$ are obtained by the average of values got in the $43$ runs in which $WS$ was able to restore the reference assignments, the analogous percentage for $MSCH$ are obtained by the average of 7 runs only.

\begin{table}[!htb]
\centering
\caption{Number of times each model was able to restore the considered reference assignments. The values marked by * are computed out of 28 possible runs, while all the others are computed out of 50 runs\label{CompatibleModels}}
\resizebox{0.5\textwidth}{!}{
\begin{tabular}{ccccc}
		\hline
		\hline
    $\#$ reference assignments  & {WS} & {CH} & {MSCH} & {MLR} \\
    \hline
    \hline
    {14} & 43 & 50 & 7  & 50 \\
    {16} & 31 & 50 & 19 & 50 \\
    {18} & 20 & 50 & 30 & 50 \\
    {20} & 15 & 50 & 35 & 50 \\
    {22} & 8  & 50 & 42 & 50 \\
    {24} & 8  & 50 & 42 & 50 \\
    {26} & 2  & 50 & 48 & 50 \\
    {27} & 1*  & 28* & 27* & 28* \\
    \hline
	\hline
    \end{tabular}%
}
\end{table}

\begin{table}[!htb]
\centering
\caption{Average percentage of countries assigned to the correct class considering the four methods at hand\label{ComparisonTable}}
\resizebox{1\textwidth}{!}{
\begin{tabular}{ccccc}
		\hline
		\hline
    $\#$ reference assignments  & {WS ($p=50$/$p=75$/$p=100$)} & {CH ($p=50$/$p=75$/$p=100$)} & {MSCH ($p=50$/$p=75$/$p=100$)} & {MLR} \\
    \hline
		\hline
    {14} & 54.8 / 48.8 / 39.0 & 57.5 / 52.0 / 43.0 & 60.2 / 49.0 / 38.8 & 40.1 \\
    {16} & 54.5 / 50.0 / 43.3 & 63.1 / 59.3 / 50.3 & 64.9 / 54.4 / 43.0 & 45.0 \\
    {18} & 57.5 / 53.0 / 49.0 & 65.2 / 61.6 / 56.8 & 68.6 / 62.3 / 52.7 & 44.6 \\
    {20} & 63.3 / 60.0 / 54.2 & 66.5 / 64.3 / 58.3 & 70.3 / 63.6 / 57.1 & 46.7 \\
    {22} & 54.1 / 52.1 / 52.1 & 66.0 / 65.3 / 63.0 & 70.2 / 62.7 / 56.4 & 43.0 \\
    {24} & 50.0 / 50.0 / 46.9 & 53.5 / 51.5 / 50.0 & 58.9 / 47.6 / 40.5 & 38.0 \\
    {26} & 50.0 / 50.0 / 50.0 & 63.0 / 63.0 / 63.0 & 66.6 / 65.6 / 59.4 & 41.0 \\
    {27} & 0.000 / 0.000 / 0.000 & 60.7 / 60.7 / 60.7 & 62.9 / 63.0 / 50.3 & 39.3 \\
		\hline
		\hline
    \end{tabular}%
}
\end{table}

Looking at the data in Fig. \ref{ComparisonFigure}, one can obviously conclude that $MLR$ is the worst in terms of predictive capability among the four considered models since, apart from very few cases, the other three methods report better percentages of correct assignments. Looking at the three methods different from $MLR$, it seems that the $WS$ is the worst since at first, as underlined above, it is not always able to represent the reference assignments and, secondly, in the cases in which it is able to represent the reference assignments the percentage of correct assignments is almost always lower than the corresponding percentage for $CH$ and $MSCH$ being able to take into account the interaction between pairs of elementary criteria. \\
Finally, comparing $CH$ and $MSCH$, it seems that the second presents better percentages of correct assignments in the case $p=50$, they are comparable for $p=75$ and $CH$ seems better than $MSCH$ for $p=100$. 


\section{Conclusions}\label{concl}%
Several methods have been proposed in Multiple Criteria Decision Aiding (MCDA) to deal with sorting problems in which alternatives need to be assigned to one among the considered classes ordered from the worst to the best with respect to the preferences provided by the Decision Maker (DM). Despite many real world sorting problems involve alternatives evaluated on interacting criteria, very few works have been developed in recent years to cope with this type of problems. In this paper, we propose a sorting method based on the Choquet integral preference model able to take into account the possible negative and positive interactions between the criteria at hand. A procedure to look for the minimal number of pairs of interacting criteria is also presented here for the first time and we think that this is an important introduction not only for our method but also for all the other methods having the Choquet integral as underlying preference model.\\
 To deal also with problems in which criteria are not interacting but also structured in a hierarchical way, the new proposal puts the Multiple Criteria Hierarchy Process (MCHP) and the Choquet integral preference model under a unified framework. The application of the MCHP to this context permits to get finer assignments recommendations not only at comprehensive level (considering all criteria simultaneously) but also a particular aspect being the most interesting for the DM. 

The proposed method permits the DM to provide preference information in terms of reference assignments, pairwise comparisons of alternatives and comparisons between criteria in terms of their importance as well as in terms of possible interactions. In this case, more than one model could be compatible with such information and, therefore, to take into account all the compatible models providing, therefore, more robust conclusions, the Robust Ordinal Regression (ROR) and the Stochastic Multicriteria Acceptability Analysis (SMAA) are applied to the new method. On one hand, ROR applied to the new sorting method provides necessary and possible assignments at comprehensive level as well as on a particular macro-criterion. On the other hand, based on a sampling procedure, SMAA gives the Class Acceptability Index (CAI) being the frequency with which an alternative is assigned to a particular class. To aggregate the different CAIs, a procedure based on the minimization of a misclassification is proposed. It permits to summarize the different assignments into a unique assignment being representative of the different CAIs. 

The merits of the proposed hierarchical, interacting and robust sorting method have been shown by its application to the financial rating of 28 European Countries. In particular, we have shown that a weighted sum is not able to represent the considered assignments, while our model is able to represent them. Moreover, we compared the predictive capabilities of two versions of our method (the one considering a minimal set of pairs of interactions and the one admitting all possible interactions between pairs of criteria) with other two methods based on a weighted sum and on a multiple linear regression model, respectively. The two versions of our method present better results than the other two since in all cases they are able to represent the reference assignments provided by the DM and, moreover, they can assign the non reference Countries to the correct class more frequently than the other two methods.  

To conclude, we think that the proposed method can be successfully applied to any MCDA sorting problem since it is able to take into account their main characteristics that is, interactions between criteria, hierarchy of criteria and robustness concerns. As further directions of research, we plan to extend the new sorting method to take into account also non-monotonic criteria and to apply it to some real world sorting problem. 

\section*{Acknowledgments}
The authors wish to acknowledge the funding by the research project ``Data analytics for entrepreneurial ecosystms, sustainable development and well being indices'' of the Department of Economics and Business of the University of Catania. 

\section*{References}

\bibliographystyle{plainnat}
\bibliography{Full_bibliography,Full_bibliographySort}


\end{document}